\documentclass{amsart}

\usepackage[utf8]{inputenc}
\usepackage[detect-all]{siunitx}

\usepackage{graphicx}
\usepackage{amsmath}
\usepackage{amsfonts}
\usepackage{amssymb}
\usepackage{a4wide}
\usepackage{listings}

\newcommand{\fv}{\vec f}
\newcommand{\uv}{\vec u}
\newcommand{\vv}{\vec v}
\newcommand{\nv}{\vec n}
\newcommand{\tv}{\vec t}
\newcommand{\Ev}{\vec E}
\newcommand{\sigmat}{\underline{\sigma}}
\newcommand{\sigmav}{\vec{\sigma}}
\newcommand{\tauv}{\vec{\tau}}
\newcommand{\taut}{\underline{\tau}}
\newcommand{\eps}{\varepsilon}
\newcommand{\epst}{\underline{\varepsilon}}
\newcommand{\epsilont}{\underline{\epsilon}}
\newcommand{\St}{\underline{S}}
\newcommand{\Ct}{\underline{C}}
\newcommand{\Dv}{\vec D}
\newcommand{\dt}{\underline{d}}
\newcommand{\et}{\underline{e}}
\newcommand{\gt}{\underline{g}}

\newcommand{\spaceV}{V}
\newcommand{\spaceW}{W}
\newcommand{\spaceSigma}{\Sigma}
\newcommand{\spaceD}{D}
\newcommand{\spaceX}{X}
\newcommand{\spacePhi}{\Phi}

\newcommand{\opdiv}{\operatorname{div}}

\newcommand{\rev}[1]{{#1}}
\newcommand{\oldnew}[2]{{#2}}
\author{Astrid S.\ Pechstein, Martin Meindlhumer and Alexander Humer}

\title{New mixed finite elements for the discretization of piezoelectric structures or macro-fibre composites}

\email{astrid.pechstein@jku.at}


\begin{document}

\begin{abstract}
    We propose a new three-dimensional formulation based on the mixed Tangential-Displacement Normal-Normal-Stress (TDNNS) method for elasticity. In elastic TDNNS elements, the tangential component of the displacement field and the normal component of the stress vector are degrees of freedom and continuous across inter-element interfaces. TDNNS finite elements have been shown to be locking-free with respect to shear locking in thin elements, which makes them suitable for the discretization of laminates or macro-fibre composites. In the current paper, we extend the formulation to piezoelectric materials by adding the electric potential as degree of freedom. 
\end{abstract}

\maketitle

\section{Introduction}

The simulation of smart, piezoelectric structures is of high interest in science and applications. A powerful method for the approximate solution  of the underlying coupled electro-mechanical equations is the finite element (FE) method. 

First FE simulations of piezoelectric structures were carried out by Allik and Hughes \cite{AllikHughes:1970} and later by Lerch \cite{Lerch:1988,Lerch:1990}. They provide volume  elements based on the principle of virtual works where the mechanical displacements and the electric potential are chosen as degrees of freedom. These finite element methods are very flexible -- in principle, they can be used to model almost any technical application. A severe drawback is the complexity of the underlying numerical system. Due to locking, flat layered piezoelectric structures have to be resolved by a sufficient number of well-shaped elements. This easily leads to computational systems with millions of unknowns even for simple applications as, e.g., thin piezoelectric patches.

Two different ways to circumvent this problem are pursued nowadays: the design of locking-free volume elements and the derivation of equations for layered plates, beams and shells. For both categories, we distinguish methods based on the principle of virtual works, and so-called mixed methods based on Hellinger-Reissner type formulations. In the former class of methods, the displacement field and the electric potential are considered as unknowns, while, in the latter class, the mechanical stresses and sometimes also the dielectric displacement field, are added as degrees of freedom.

For both volume elements as well as layered plate, beam or shell elements, it has been shown numerically that mixed methods provide good results for coarse discretizations  independently of the layer thickness. We mention the volume element by \oldnew{Sze and Yao \cite{SzeYao:2000}}{Sze, Yao and Yi \cite{SzeYaoYi:2000}}, the solid shell element by Klinkel and Wagner \cite{KlinkelWagner:2006} and the geometrically nonlinear element by Ortigosa and Gil \cite{OrtigosaGil:2016}. Reissner-type mixed zigzag formulations were successfully used by \cite{CarreraBoscolo:2007,CarreraNali:2009,WuLin:2014}. 

These findings motivate our suggestion for a new family of piezoelectric elements. In \cite{PechsteinSchoeberl:11,PechsteinSchoeberl:12} the ``Tangential Displacement Normal Normal Stress'' (TDNNS) finite element method was introduced for linear elastic solids. The elements are based on a Hellinger-Reissner formulation, where the tangential component of the displacements as well as the normal component of the (normal) stress vector are considered as degrees of freedom. In \cite{PechsteinSchoeberl:12} it was shown that these elements are locking-free when used as flat prismatic elements. In the current contribution, we propose an extension of these elements to piezoelectric materials.


We discuss the implementation of the proposed elements in the open-source software package Netgen/NGSolve for the case of a bimorph beam.
We show the accuracy and convergence rates for this exemplary problem. We present results for the more advanced problem of computing effective material properties of a  $d_{15}$ MFC.


\section{The problem of linear piezoelasticity}

Let $\Omega \subset \mathbb{R}^3$ describe a solid, which is made of elastic, piezoelectric material. In the following, we derive a formulation for linear piezoelasticity, i.e. we assume linearity of the (piezo-)elastic material laws as well as the case of small deformations. 
This simplest form of electro-mechanical coupling, which describes the behavior of piezoelectric materials for a given poling state, is also referred to as ``Voigt's linear theory of piezoelectricity"~\cite{Kamlah:2001}.
Further, we neglect the electrically induced contributions to the mechanical balance laws, which preserves the symmetry of the Cauchy stress tensor $\sigmat$.  
%
%
We are interested in finding the displacement field $\uv$ and the electric potential $\phi$ subject to body forces $\fv$ and (suitable) boundary conditions.
Derived from these fields are the electric field $\Ev = - \nabla \phi$ and the linear strain tensor $\epst = \tfrac12 (\nabla \uv + \nabla \uv^T)$. 
Standard finite element formulations are based on a variational principle such as the principle of virtual works or D'Alembert's principle. In such a formulation, displacements $\uv$ and electric potential $\phi$ are considered as independent variables. The degrees of freedom of the elements represent these fields. We call $\uv$ and $\phi$ the \emph{primal quantities}.

Opposed to the primal quantities, the \emph{dual quantites} of interest are the Cauchy stress tensor $\sigmat$ and the dielectric displacements $\Dv$. Usually, the quantities $\sigmat$ and $\Dv$ are computed in a postprocessing step from $\uv$ and $\phi$, using the constitutive laws. This implies that the order of approximation for $\sigmat$ and $\Dv$ is one less than for $\uv$ and $\phi$.

\begin{figure}
\begin{center}
\includegraphics[width=0.8\textwidth]{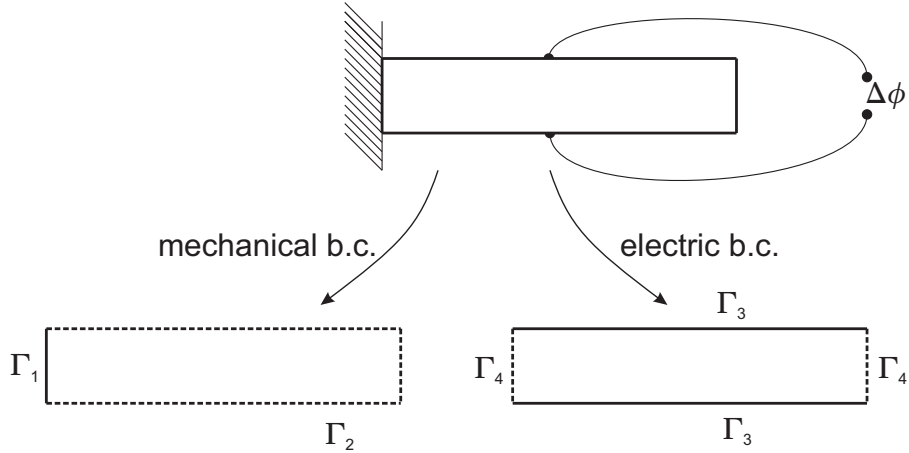}
\end{center}
\caption{Visualization of mechanic and electric boundary conditions on a clamped piezo beam which is electroded on top and bottom.}\label{fig:bc}
\end{figure}

Both dual quantities satisfy a balance equation of divergence form: for the stress tensor, we have the mechanical balance equation. The dielectric displacements satisfy Gauss' law. 
\begin{align}
-\opdiv \sigmat &=\ \fv && \text{in } \Omega, \label{eq:balance1}\\
-\opdiv \Dv &=\ \rho_e &&\text{in } \Omega. \label{eq:balance2}
\end{align}
In the above relation, $\rho_e$ denotes the free charge density, which vanishes ($\rho_e = 0$) for non-conducting solids as, e.g., piezoelectric ceramics.
The mechanical boundary conditions are 
\begin{align}
\uv &= \vec 0 \text{ on } \Gamma_1 & \text{and} && \sigmav_n &= \tv_n \text{ on } \Gamma_2 = \partial \Omega \backslash \Gamma_1.
\end{align}
The electrical boundary conditions are
\begin{align}
\phi &= \phi_0\text{ on } \Gamma_3 & \text{and} && D_n &= q_0 \text{ on } \Gamma_4 = \partial \Omega \backslash \Gamma_3.
\end{align}
A visualization of boundary conditions for a simple example of a clamped piezo beam can be found in Figure~\ref{fig:bc}.

\subsection{Different formulations of the constitutive laws}
Stress and strain are second order symmetric tensors, which are represented by symmetric three-by-three matrices.
In the following, we use Voigt's notation for stresses and strains, where $\epst$ and $\sigmat$ are interpreted as six-dimensional vectors. In our notation, we do not distinguish between symmetric matrix and vector, as it will be clear from the context which one is to use.

There are several ways to formulate the material laws of linear piezoelasticity, which are equivalent for linear materials. For standard finite element formulations, one usually has the dual quantities $\sigmat$ and $\Dv$ depending on the primal quantities $\epst$ and $\Ev$. This results in 
\begin{align}
\sigmat &=\ \Ct^E \epst - \et \Ev, \label{eq:etype1}\\
\Dv &= \ \et^T \epst + \epsilont^{\eps} \Ev. \label{eq:etype2}
\end{align}
Here, $\Ct^E$ denotes the elasticity tensor measured at constant electric field and $\epsilont^{\eps}$ is the dielectric tensor or electric permittivity at constant mechanical strain. The piezoelectric coupling is described by the piezoelectric permittivity tensor $\et$. 

Another widely used set of material parameters uses the piezoelectric tensor $\dt$. Then, one additionally needs the compliance $\St^E$ and the dielectric tensor at constant mechanical stresses $\epsilont^{\sigma}$,
\begin{align}
\epst &=\ \St^E \sigmat + \dt^T \Ev, \label{eq:dtype1}\\
\Dv &=\ \dt \sigmat + \epsilont^{\sigma} \Ev. \label{eq:dtype2}
\end{align}

Of course, the material parameters are connected by the well-known relations
\begin{align}
\St^{E} &=\ (\Ct^{E})^{-1},& \dt &=\ \et \St^{E}, & \epsilont^{\sigma} &=\ \epsilont^{\eps} + \dt \et^T.
\end{align}

Less often, one finds the piezoelectric tensor $\gt$. Using the latter, strain and electric field can be expressed depending on stresses and dielectric displacements,
\begin{align}
\epst &=\ \St^D \sigmat + \gt^T \Dv,\label{eq:gtype1}\\
\Ev &=\ -\gt \sigmat + (\epsilont^{\sigma})^{-1} \Dv. \label{eq:gtype2}
\end{align}
Here we use
\begin{align}
\gt &=\ (\epsilont^{\sigma})^{-1} \dt, & \St^{D} &=\ \St^{E} - \dt^T\gt.
\end{align}
We will use the $\dt$-type and the $\gt$-type formulations for the proposed mixed finite elements\oldnew{.}{, as  then strain and, in the second variant, also electric field are provided as functions of the dual quantities stress and dielectric displacement. Note that, for linear materials, the $\dt$-type and $\gt$-type formulations are equivalent to the $\et$-type formulation, and thus always available. 
}

\section{Preliminaries for the mixed finite element method}

To develop a finite element method, we assume $\mathcal{T} = \{T\}$ to be a finite element mesh of the domain $\Omega$, consisting of tetrahedral, prismatic or hexahedral elements. Of course, all results of this contribution can be transferred to two-dimensional problems using triangular or quadrilateral elements. By $\nv$ we denote the outward unit normal on the (element or domain) boundary $\partial T$ or $\partial \Omega$. 
On each element or domain boundary surface, a general vector field $\vv$ can be split into normal and tangential components by $\vv = v_n \nv + \vv_t$ with $v_n = \vv\cdot \nv$ and $\vv_t = \vv - v_n \nv.$ Note that the normal component $v_n$ is scalar, while the tangential component $\vv_t$ is vectorial. Any tensor field $\taut$ has a normal vector $\tauv_n = \taut \nv$ on a surface, which can again be split into normal and tangential components $\tau_{nn} = \tauv_n\cdot n$ and $\tauv_{nt} = \tauv_n - \tau_{nn} \nv$.


The TDNNS method is a mixed finite element method\oldnew{ based on Reissner's principle}{, which can be seen as a variant of Reissner's principle} \rev{\cite{Reissner:1950}}. Displacements and stresses are considered as unknowns, see \cite{PechsteinSchoeberl:11,PechsteinSchoeberl:16}. 
The tangential component of the displacement $\uv_t$ and the normal component of the stress vector $\sigma_{nn}$ are degrees of freedom of the finite element. 
\rev{These quantities are also the essential boundary conditions of the finite element method, the finite element functions explicitely satisfy
\begin{align}
\uv_t &= 0 \text{ on } \Gamma_1, & \sigma_{nn} &= t_{nn} \text{ on } \Gamma_2.
\end{align}
On the other hand, natural boundary conditions on $u_n$ and $\sigmav_{nt}$ will enter, if non-zero, the right hand side of the variational formulation as external works.}
We use finite element spaces for which \oldnew{these respective quantities}{these degrees of freedom} are continuous across element interfaces. Note that, for this choice of degrees of freedom, the finite element displacement field can be discontinuous. A gap between elements in normal direction may open up, while sliding in the tangential direction is prohibited. In the solution, gaps are controlled by an extra interface term in the principle of virtual works, see \eqref{eq:TDNNS} and \eqref{eq:defdiv1}--\eqref{eq:defdiv2}.

N\'ed\'elec \cite{Nedelec:80,Nedelec:86} introduced tangential-continuous finite elements, which are commonly used to describe the electric field in Maxwell's equations. We use the elements from \cite{Nedelec:86} for the displacements. For the stresses, we introduced normal-normal continuous elements in \cite{PechsteinSchoeberl:11}. On simplicial meshes, these finite element spaces can be shortly described by
\begin{eqnarray}
\uv, \delta \uv \in \spaceV_h &=& \{ \vv: \vv|_T \in [P^k(T)]^3, \vv_t \text{ continuous}\},  \label{eq:spaceVh}\\
\sigmat, \delta \sigmat \in \spaceSigma_h &=& \{ \taut: \taut|_T \in [P^k(T)]^{3\times3}_{sym}, \tau_{nn} \text{ continuous} \}.  \label{eq:spaceSigmah}
\end{eqnarray}
Here, $P^k(T)$ denotes the space of polynomials of order at most $k$ on (simplicial) element $T$. Where prismatic elements are concerned, these spaces are extended exploiting their tensor product nature.
For the tangential continuous space $\spaceV_h$, the corresponding elements and shape functions we use are described in \cite[p.\ 92f.]{Zaglmayr:2006}.
For the normal-normal continuous stress space, tetrahedral, prismatic and hexahedral elements are provided in \cite{PechsteinSchoeberl:11}, additionally two-dimensional triangular or quadrilateral elements exist.  In \cite{PechsteinSchoeberl:12} it was shown that the method works well for thin prismatic or hexahedral elements. All these elements are implemented in the open-source software package Netgen/NGSolve\footnote{Open-source software package Netgen/NGSolve \tt https://ngsolve.org}. \rev{In Figure~\ref{fig:dofs}, we illustrate the interface degrees of freedom of the displacement and stress elements of polynomial order one. Summing up, for the stresses, we have 24 dof, for the displacements 30 coupling dof. Note that all internal degrees of freedom can be eliminated while assembling the finite element matrix. Also, though the elements sport more degrees of freedom than classical nodal elements, the coupling through element edges and faces is much weaker than nodal coupling. Thus the stiffness matrix is sparser, and can be solved faster by a direct solver.}

\rev{
\begin{figure}
\begin{center}
\includegraphics[width=0.6\textwidth]{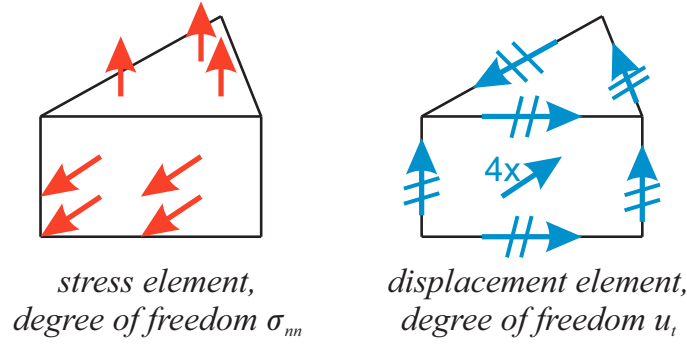}
\caption{Illustration of the prismatic TDNNS stress and displacement elements of polynomial order one.}\label{fig:dofs}
\end{center}
\end{figure}
}

We provide a variational formulation for the TDNNS method, skipping the details which can be found in \cite{PechsteinSchoeberl:11, PechsteinSchoeberl:12, PechsteinSchoeberl:16}.
The formulation is based on Reissner's principle, and reads: find $\uv \in \spaceV_h$ satisfying $\uv_t = 0$ on $\Gamma_1$ and $\sigmat \in \spaceSigma_h$ satisfying $\sigma_{nn} = t_{nn}$ on $\Gamma_2$ such that
\begin{align}
\int_\Omega \St \sigmat : \delta \sigmat\, d\Omega - \langle \epst(\uv), \delta \sigmat\rangle -\langle \epst(\delta \uv), \sigmat \rangle &=\ \rev{-}\int_\Omega \fv \cdot \delta \uv\,d\Omega\ \oldnew{+}{ -} \int_{\Gamma_2} \tv_{nt} \delta \uv_t\, d\Gamma ,
\label{eq:TDNNS}
\end{align}
for all virtual displacements $\delta \uv \in \spaceV_h$ and virtual stresses $\delta \sigmat \in \spaceSigma_h$ which satisfy the corresponding homogeneous essential boundary conditions. We note that the tangential displacement $\uv_t$ and normal stress $\sigma_{nn}$ are the essential degrees of freedom. Normal displacement $u_n$ and shear stress $\sigmav_{nt}$ are natural boundary conditions. Inhomogeneous conditions on the shear stress $\sigmav_{nt} = \tv_{nt}$ and -- if applicable -- also of the normal displacement are added to the right hand side of \eqref{eq:TDNNS}.


In the variational principle \eqref{eq:TDNNS}, we see duality products of the form $\langle \epst(\uv), \sigmat\rangle$ instead of integrals of the form $\int_\Omega \epst(\uv) : \sigmat\, d\Omega$ in common methods. This distinction is necessary to be mathematically correct, since the strain $\epst(\uv)$ of a finite element function $\uv$ is a distribution. Recall that a finite element displacement function $\uv$ is not completely continuous, but gaps in the normal displacement may arise. These gaps lead to an additional distributional part of the strain, which is  evident as element-wise surface integrals in formulas \eqref{eq:defdiv1}--\eqref{eq:defdiv2}. The duality product is well defined only if the stress field $\sigmat$ is normal-normal continuous. Note that this is exactly the defining property of the stress elements, see \eqref{eq:spaceSigmah}. For finite element functions on the mesh $\mathcal{T}$, the duality product can be evaluated element-wise by volume and surface integrals. The surface integrals represent the distributional terms on element interfaces mentioned above. The following two formulas are equivalent, and motivate that the duality product can also be viewed as the (negative) distributional divergence of the stress tensor,
\begin{eqnarray}
\langle\epst(\uv), \sigmat \rangle &=& \sum_{T \in \mathcal{T}} \Big( \int_T \sigmat : \epst(\uv)\, d\Omega - \int_{\partial T} \sigmat_{nn} \uv_n\, d\Gamma \Big) \label{eq:defdiv1}\\
&=& \sum_{T \in \mathcal{T}} \Big( -\int_T \opdiv \sigmat \cdot \uv\, d\Omega + \int_{\partial T} \sigmat_{nt} \uv_t\, d\Gamma \Big)\ = \
-\langle \opdiv \sigmat, \uv\rangle. \label{eq:defdiv2}
\end{eqnarray}
The equivalence of formulas \eqref{eq:defdiv1} and \eqref{eq:defdiv2} can be shown by integration by parts on each element, and using the continuity of $\sigma_{nn}$ and $\uv_t$, respectively. For a more involved mathematical motivation see \cite{PechsteinSchoeberl:11,PechsteinSchoeberl:16}.

\section{Mixed finite elements for piezoelastic structures}

In the sequel, we present two piezoelectric finite elements based on the elastic TDNNS method. In the first variant, the electric field $\Ev = - \nabla \phi$ is added by considering the electric potential $\phi$ as a further unknown. Additionally, we propose a method where the electric potential $\phi$ and the dielectric displacements $\Dv$ are added as independent variables. 
Numerical results for both methods shall be presented, indicating that while the former method is probably easier to implement, the latter yields more accurate results. 

\subsection{Revisiting standard piezoelectric elements}

We shortly discuss the standard variational formulation based on the principle of virtual works. In this formulation, the displacements $\uv$ and the electric potential $\phi$ are considered independent unknowns. One plugs the $\et$-type material laws \eqref{eq:etype1}--\eqref{eq:etype2} into the balance equations \eqref{eq:balance1}--\eqref{eq:balance2},
\begin{align}
-\opdiv( \Ct^{E} \epst + \et^T \nabla\phi) &= -f,\label{eq:primal1}\\
-\opdiv(\et \epst - \epsilont^{\eps} \nabla \phi) &= 0. \label{eq:primal2}
\end{align}
To derive a variational formulation, one multiplies the first equation \eqref{eq:primal1} by a virtual displacement $\delta \uv$ satisfying the (zero) boundary condition on $\Gamma_1$, and the second equation \eqref{eq:primal2} by a virtual potential satisfying the (zero) boundary condition on $\Gamma_3$. Then one integrates by parts and employes the natural boundary conditions on $\Gamma_2$ and $\Gamma_4$, leading to
\begin{align}
\int_\Omega (\Ct^{E} \eps + \et^T \nabla \phi) : \delta \eps\, d\Omega + \int_\Omega (\et \eps - \epsilont^{\eps} \nabla \phi)& \cdot \delta \nabla \phi\, d\Omega =
\int_\Omega \sigma : \delta \epst\,d\Omega - \int_\Omega \Dv \cdot \delta\Ev\, d\Omega\\&=\ \int_\Omega \fv \cdot \delta \uv\, d\Omega +  \int_{\Gamma_2} \tv_{nt} \cdot \delta \uv_t\, d\Gamma + \int_{\Gamma_4} q_0 \delta \phi\, d\Gamma.
\end{align} 

We use standard continuous (e.g.\ nodal) finite elements for the displacements and the electric potential, which satisfy the respective (homogenized) boundary conditions,
\begin{align}
\uv, \delta \uv &\in\ \{ \uv_h \in [\spaceW_h]^3: \uv_h = 0 \text{ on } \Gamma_1\},\\
\phi, \delta \phi &\in\ \{ \phi_h \in \spaceW_h: \phi_h = 0 \text{ on } \Gamma_3\},\\
\text{with }  \spaceW_h &:=\ \{ w: w|_T \in P^k(T), w \text{ continuous}\}. \label{eq:wh}
\end{align}
Again, the definition of $\spaceW_h$ given in \eqref{eq:wh} holds for simplicial elements. It can be extended to prismatic tensor-product elements. In our computations, we use the shape functions described in \cite[p.\ 95f.]{Zaglmayr:2006}, which are implemented in Netgen/NGSolve.

\subsection{TDNNS-based elements using the electric potential} 

We shall now develop a first mixed formulation for piezoelasticity, which is based on the TDNNS formulation. In the numerical examples, the formulation of this section is indicated as ``first variant'' or V1. The independent unknowns are the displacement vector $\uv \in \spaceV_h$, the stess tensor $\sigmat \in \spaceSigma_h$ and the electric potential $\phi \in \{ \phi_h \in \spaceW_h: \phi_h = 0 \text{ on } \Gamma_3\}$. The essential degrees of freedom are the tangential displacement, the normal component of the stress vector, and the nodal values of the electric potential.

We use the $\dt$-type material laws \eqref{eq:dtype1}--\eqref{eq:dtype2}, and eliminate the dielectric displacements from the balance equations \eqref{eq:balance1}--\eqref{eq:balance2},
\begin{align}
- \St^{E} \sigmat + \dt^T \nabla \phi + \epst &= 0,\\
- \opdiv \sigmat &= f,\\
-\opdiv(\dt \sigmat - \epsilont^{\sigma} \nabla \phi) &= 0. \label{eq:lastline}
\end{align}
Then we multiply the first line by a virtual stress $\delta \sigma$, the second line by a virtual displacement $\delta \uv$ and the third line by a virtual potential $\delta\phi$ with $\delta \phi = 0$ on $\Gamma_3$.
We use the distributional strain and divergence operators for the mechanical quantities \eqref{eq:defdiv1}--\eqref{eq:defdiv2}. 
\begin{align}
\int_\Omega (-\St^{E} \sigmat + \dt^T \nabla \phi): \delta \sigmat\,d\Omega  + \langle\epst, \delta \sigmat\rangle &= 0,\\
\langle \delta \epst, \sigmat\rangle &= \int_\Omega \fv \cdot \delta \uv\, d\Omega + \int_{\Gamma_2} \tv_{nt} \cdot \delta \uv_t\, d\Gamma.
\end{align}
In the last line \eqref{eq:lastline}, we apply integration by parts in the same way as for the standard elements. 
\begin{align}
-\int_\Omega \opdiv(\dt \sigmat - \epsilont^{\sigma} \nabla \phi)\delta \phi\,d\Omega &=
\int_\Omega (\dt \sigmat - \epsilont^{\sigma} \nabla \phi) \cdot \nabla \phi\, d\Omega - \int_{\Gamma_4} q_0\, \delta \phi\, d\Gamma = 0.
\end{align}
Here we used that $\delta \phi$ vanishes on $\Gamma_3$ and $D_n = q_0$ on the remainder $\Gamma_4$. 
Summing up, we arrive at
\begin{align}
-\int_\Omega (\St^{E} \sigmat - \dt^T\nabla \phi) : \delta \sigmat\, d\Omega 
+ \langle \epst(\uv), \delta \sigmat\rangle + \langle \epst(\delta \uv), \sigmat \rangle +
\int_\Omega (\dt \sigmat - \epsilont^{\sigma} \nabla \phi) \cdot \delta\nabla \phi\, d\Omega =&\\
-\int_\Omega \epst(\sigmat,\Ev) : \delta \sigmat\, d\Omega + \langle \epst(\uv), \delta \sigmat\rangle + \langle \epst(\delta \uv), \sigmat \rangle - \int_\Omega \Dv\cdot \delta \Ev\, d\Omega
 =&\\
= \int_\Omega \fv \cdot \delta \uv\,d\Omega + \int_{\Gamma_2} \tv_{nt} \cdot \delta \uv_t\, d\Gamma +\int_{\Gamma_4} q_0\, \delta \phi\, d\Gamma. \label{eq:V1}
\end{align}

The performance of thin prismatic elements is tested in the sequel. We shall see that it is free from locking if the order of the electric potential is at least $k=2$, while we can choose linear elements for the mechanical quantities. For the lowest order case, the quadratic behavior of the electric potential in thickness direction cannot be represented by the degrees of freedom, as we have only a linear behavior of $\phi$. Therefore, we expect a deterioration of accuracy when using only one layer of elements in thickness direction. However, for higher orders, the quadratic variation can be represented, and accurate results are obtained.

\subsection{TDNNS-based elements using the electric potential and the dielectric displacements}

We propose a variant of the finite element method, which involves both dual quantities mechanical stresses and dielectric displacements as unknowns. Thus, we end up with four independent unknown fields displacement $\uv$, electric potential $\phi$, stress $\sigmat$ and dielectric displacement $\Dv$. We will first derive the variational equations. From this derivations, we will deduce the appropriate degrees of freedom for the electric quantities, as well as the essential boundary conditions. In the numerical results, this finite element formulation will be referred to as ``second variant'' or V2.

We use the material laws in $\gt$-type form \eqref{eq:gtype1}--\eqref{eq:gtype2}, and both balance equations,
\begin{align}
-\St^{D} \sigmat - \gt^T \Dv + \epst &=\ 0,\\
-\gt \sigmat + (\epsilont^{\sigma})^{-1} \Dv + \nabla \phi &=\ 0,\\
-\opdiv \sigmat &= \fv,\\
-\opdiv \Dv &=\ 0.
\end{align}
We multiply the first line by a virtual stress $\delta \sigmat \in \spaceSigma_h$ and the third line by a virtual displacement $\delta \uv \in \spaceV_h$, which satisfy the corresponding homogeneous boundary conditions for $\uv_t$ and $\sigma_{nn}$.  The second line is multiplied by a virtual dielectric displacement $\delta \Dv$ which satisfies $\delta D_n = 0$ on the insulated boundary part $\Gamma_4$. The last line is multiplied by a virtual potential $\delta \phi$ which satisfies no boundary condition a priori. Integrating over the domain and using the distributional strain and divergence operators for the mechanical quantities, we arrive at
\begin{align}
-\int_\Omega(\St^{D} \sigmat - \gt^T \Dv):\delta \sigmat\, d\Omega + \langle\epst(\uv),\delta \sigma\rangle &=\ 0,\\
-\int_\Omega(\gt \sigmat + (\epsilont^{\sigma})^{-1} \Dv)\cdot \delta \Dv\, d\Omega + \int_\Omega \nabla \phi\cdot \delta \Dv\, d\Omega &=\ 0, \label{eq:tdnnsdualherleitung2}\\
\langle\delta \epsilont, \sigmat\rangle &= \int_\Omega \fv \cdot \delta \uv\, d\Omega + \int_{\Gamma_2} \tv_{nt} \cdot \delta \uv_t\, d\Gamma,\\
-\int_\Omega\opdiv \Dv\, \delta \phi\, d\Omega &=\ 0.
\end{align}

In the next step, we apply standard integration by parts in the second integral of eq.\ \eqref{eq:tdnnsdualherleitung2}. We use the boundary conditions $\delta D_n = 0$ on $\Gamma_4$ and $\phi = \phi_0$ on $\Gamma_3$. Moreover, we assume that $\delta \Dv$ is smooth enough such that it allows for a divergence, i.e.\ $\opdiv \Dv$ exists at least in $L^2$ sense. We will comment on this condition below, as it motivates the choice of finite element degrees of freedom.
In this case, we have
\begin{align}
\int_\Omega \nabla \phi\cdot \delta \Dv\, d\Omega &=\ -\int_\Omega\phi\  \opdiv\delta\Dv\, d\Omega + \int_{\Gamma_3} \underbrace{\phi}_{=\phi_0}\, \delta D_n\, d\Gamma + \int_{\Gamma_4} \phi\, \underbrace{\delta D_n}_{=0}\, d\Gamma \\
&= -\int_\Omega\phi\ \delta \opdiv\Dv\, d\Omega + \int_{\Gamma_3} \phi_0\, \delta D_n\, d\Gamma.\label{eq:divD}
\end{align}
Note that an inhomogeneous boundary condition for the electric potential $\phi = \phi_0 \neq 0$ is a natural boundary condition in this formulation, which appears at the right hand side of \eqref{eq:divD} or later \eqref{eq:divD2}. 
Inserting the identity \eqref{eq:divD} in \eqref{eq:tdnnsdualherleitung2} leads to the final variational formulation. Before we formally put it down, we discuss the finite element spaces used for dielectric displacements and electric potential.

For piecewise smooth (or polynomial) finite element functions $\Dv$, the divergence is in $L^2$ if and only if the normal component $D_n$ is continuous across element interfaces. Thus, the normal component $D_n$ has to be a degree of freedom living on element faces in 3D, or edges in 2D. Different elements satisfying this constraint were introduced. We cite the original work by Raviart and Thomas \cite{RaviartThomas:1977}, which was generalized to three dimensional problems in \cite{Nedelec:80}. For an overview on divergence-conforming elements we refer to the monograph \cite{BoffiBrezziFortin:2013}. 

In the right hand side of \eqref{eq:divD} as well as the final variational formulation, no derivatives of the electric potential occur. Thus, we use totally discontinuous elements for $\phi, \delta \phi$.
For $\Dv$, we use the divergence-conforming elements implemented in Netgen/NGSolve, which are documented in the thesis of Zaglmayr \cite{Zaglmayr:2006}. For simplicial elements, the spaces for electric potential and dielectric displacements can be described as
\begin{align}
\phi, \delta \phi \in \spacePhi^{disc}_h &= \{ \phi: \phi|_T \in P^k(T)\},\\
\Dv, \delta \Dv \in \spaceD_h &= \{ \Dv: \Dv|_T \in [P^k(T)]^3, D_n \text{ continuous} \}  \label{eq:spaceDh}.
\end{align}

As the dielectric displacement field is divergence free, the number of degrees of freedom may be reduced further. In NGSolve, there is an option to use only divergence free higher-order basis functions in the space $\spaceD_h$ above. The shape functions are then divergence free, or the divergence is constant on each element. Then, the electric potential $\phi$ can be approximated by piecewise constant finite element functions, i.e.\ we have one degree of freedom per element for $\phi$. The according spaces for simplicial elements are
\begin{align}
\phi, \delta \phi \in \spacePhi^{disc,0}_h &= \{ \phi: \phi|_T \in P^0(T)\}, \label{eq:Phi0}\\
\Dv, \delta \Dv \in \spaceD^0_h &= \{ \Dv: \Dv|_T \in [P^k(T)]^3, D_n \text{ cont.}, D_n = 0 \text{ on } \Gamma_4, \opdiv \Dv \in P^0(T) \}  \label{eq:D0}.
\end{align}
The main benefit of this option is that the element matrices and also the overall stiffness matrix is smaller and better conditioned. Degrees of freedom for the electric potential and the dielectric displacements are saved. The computed solution for the dielectric displacements is not affected by this reduction of degrees of freedom. However, the electric field cannot be evaluated via $\Ev = - \nabla \phi$, as $\phi$ is only constant per element. The material law has to be used instead,
\begin{align}
\Ev = \dt \sigmat + \epsilont^{\sigma} \Dv.
\end{align} 
Then the accuracy of the electric field is the same as that of stresses and dielectric displacements.

Using the finite element spaces above, the finite element problem is to find $\uv \in \spaceV_h$ with $\uv_t = 0$ on $\Gamma_1$, $\sigmat \in \spaceSigma_h$ with $\sigma_{nn} = t_{nn}$ on $\Gamma_2$, $\Dv \in \spaceD_h$ (or $\spaceD^{0}_h$) with $D_n = q_0$ on $\Gamma_4$ and $ \phi \in \spacePhi^{disc}_h$  (or $\spacePhi^{disc,0}_h$) such that for all virtual functions $\delta \uv \in \spaceV_h, \delta \sigmat \in \spaceSigma_h, \delta \Dv \in \spaceD_h$ (or $\spaceD^{0}_h$) and $ \delta \phi \in \spacePhi^{disc}_h$  (or $\spacePhi^{disc,0}_h$) which satisfy the respective homogeneous boundary conditions
\begin{align}
-\int_\Omega(\St^{D} \sigmat - \gt^T \Dv):\delta \sigmat\, d\Omega + \langle\epst(\uv),\delta \sigma\rangle &=\ 0, \label{eq:V21}\\
-\int_\Omega(\gt \sigmat + (\epsilont^{\sigma})^{-1} \Dv)\cdot \delta \Dv\, d\Omega - \int_\Omega \phi\cdot \delta \opdiv\Dv\, d\Omega &=\ -\int_{\Gamma_3} \phi_0\, \delta D_n\, d\Gamma, \label{eq:potbc}\\
\langle\delta \epsilont, \sigmat\rangle &= \int_\Omega \fv \cdot \delta \uv\, d\Omega + \int_{\Gamma_2} \tv_{nt} \cdot \delta \uv_t\, d\Gamma,\\
-\int_\Omega\opdiv \Dv\, \delta \phi\, d\Omega &=\ 0. \label{eq:divD2}
\end{align}

\section{Numerical results}

\subsection{Bimorph beam}
\begin{figure}
\begin{center}
\includegraphics[width=0.6\textwidth]{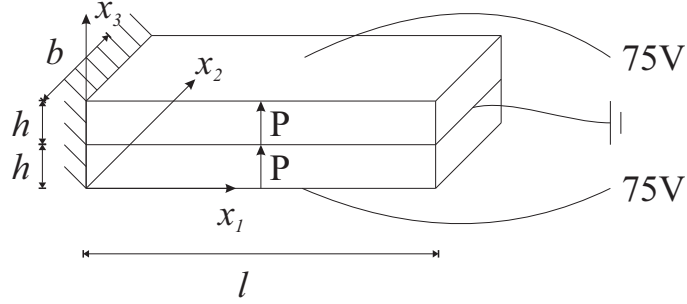}
\end{center}
\caption{Setup of the bimorph beam} \label{fig:bimbeam}
\end{figure}

The first example is a benchmark test of a piezoelectric bimorph beam. The beam is clamped at $x_1 = 0$. The length, width and height are $l = \SI{100}{\milli\meter}$, $b = \SI{10}{\milli\meter}$ and $2h = 2 \times \SI{0.5}{mm}$. The two layers of the bimorph beam are both made from PZT-5 and poled in $x_3$ thickness direction. The material parameters used are summarized in Table~\ref{tab:PZT}. The beam is electroded at the upper and lower surface, and in the interior between the layers. A constant electric potential of $\phi_0 = \SI{75}{V}$ is applied to electrodes on the upper and lower surface of the beam, while the interior electrode is grounded. See Figure~\ref{fig:bimbeam} for a sketch. 

\rev{We use both proposed variants of the method, indicating the first variant  with V1 and the second variant with V2. Recall that V1 includes only the electric potential with continuous nodal elements, while V2 uses discontinuous electric potential and normal-continuous dielectric displacements. 
We include the boundary conditions as follows: in the V1 approximation, $\uv_t = 0$ on the clamped end, $\sigma_{nn} = 0$ on all other (free) surfaces, and $\phi = \phi_0$ on the electrodes. All natural boundary conditions are homogeneous, thus no external virtual works enter the formulation. In the V2 approximation, the stress and displacement boundary conditions remain unchanged. However, now $D_n = 0$ is an essential boundary condition on all non-electroded surfaces, while $D_n$ is free on the upper and lower electrode, and free to jump across the internal electrode. The potential boundary condition now enters the right hand side as indicated in \eqref{eq:potbc}.}

\begin{table}
\caption{Material constants used for the bimorph beam.} \label{tab:PZT}
\begin{center}
\begin{tabular}{|lr|lr|}
\hline
Parameter & & Parameter &\\\hline
$C^{E}_{11}$ & $120 \times 10^9$ N/m$^2$\qquad \qquad &  $e_{31}$ & $-5.35$ C/m$^2$  \\
$C^{E}_{33}$ & $111\times 10^9$ N/m$^2$\qquad \qquad &  $e_{33}$ & $15.8$ C/m$^2$  \\
$C^{E}_{12}$ & $75.2\times 10^9$ N/m$^2$\qquad \qquad &  $e_{15}$ & $12.3$ C/m$^2$  \\
$C^{E}_{13}$ & $75.1\times 10^9$ N/m$^2$\qquad \qquad &  $\epsilon^{\eps}_{11}/\epsilon_0$ & 919  \\
$C^{E}_{44}$ & $21.1\times 10^9$ N/m$^2$\qquad \qquad &  $\epsilon^{\eps}_{33}/\epsilon_0$ & 827 \\
$C^{E}_{66}$ & $22.6\times 10^9$ N/m$^2$\qquad \qquad &   &\\\hline
\end{tabular}
\end{center}
\end{table}


%

The beam is discretized using prismatic elements in the plane. We provide a convergence study, where we provide the relative error of the average tip deflection \oldnew{(in $L^2$ sense $\|u_z - u_{z,h}\|_{L^2(\Gamma_{tip})}$)}{$|\bar u_z - \bar u^{ref}_z|/|\bar u^{ref}_z|$}, and the $L^2$ error of the displacement in the whole beam $\|\uv - \uv_h\|_{L^2(\Omega)}$. To obtain \oldnew{such an error}{the $L^2$ error}, the computed values were compared to a simulation using standard elements of higher order on the same mesh. \rev{For the tip displacement, the standard solution on the finest mesh was used for comparison. As the exact solution sports singularities at the clamped end and at the boundaries of the electrodes, an adaptive mesh refinement is necessary to get higher convergence orders for second order elements. Otherwise, the convergence is limited by the singularity, which would mean convergence of order $h^2 \simeq \text{\#ndof}$ in $L^2$ sense.} We use an adaptive mesh refinement strategy, where we employ an error estimator of Zienkiewicz-Zhu type \cite{ZienkiewiczZhu:1987}: in a postprocessing step, the (discontinuous) computed stresses are interpolated to continuous stresses. The difference between the original discontinuous and interpolated continuous stress is used as an error indicator. Elements with error indicator higher than $0.5$ times the maximum error indicator are marked for refinement. \rev{In Figure~\ref{fig:beamsigma} and Figure~\ref{fig:beamd}, we display the stress component $\sigma_{xx}$ and the dielectric displacement $D_z$ for the coarsest and most refined finite element mesh for method V2 of order $k=1$ using divergence-free high-order shape functions. One can see that the mesh is refined towards the edges of beam, where steep gradients of stress and dielectric displacement occur. As expected, the refinement towards the corners at the clamped end is strongest.

\begin{figure}
\begin{center}
\includegraphics[width=0.8\textwidth]{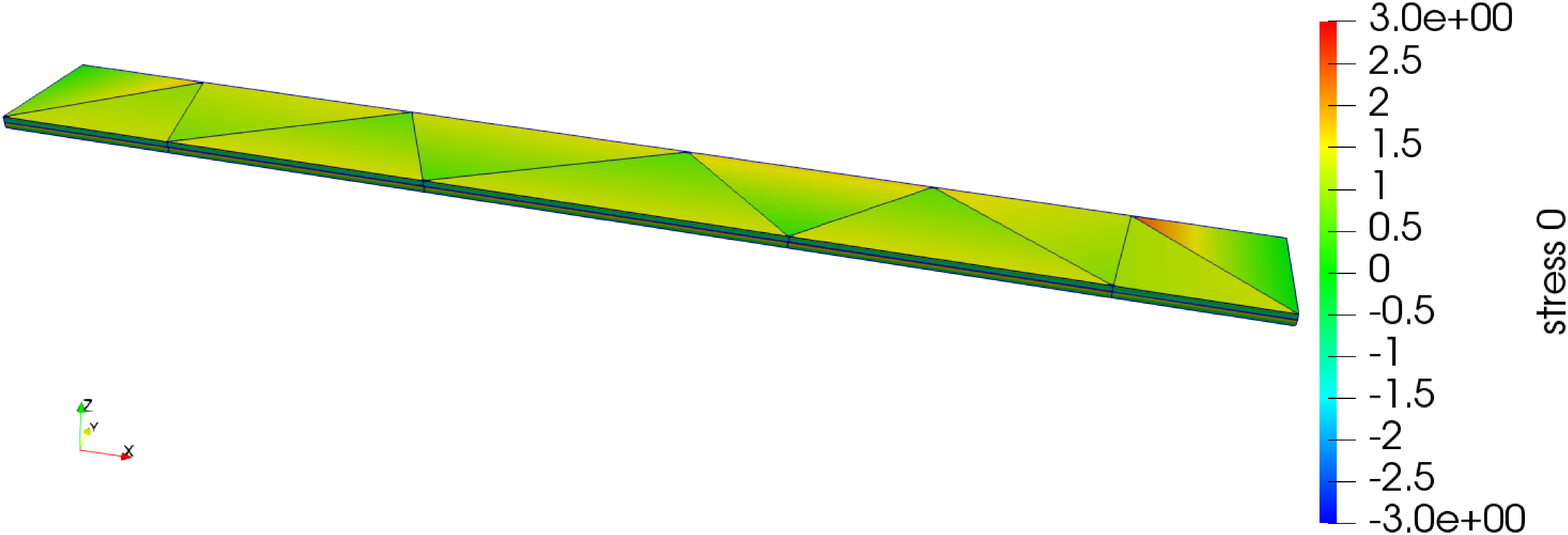}
\includegraphics[width=0.8\textwidth]{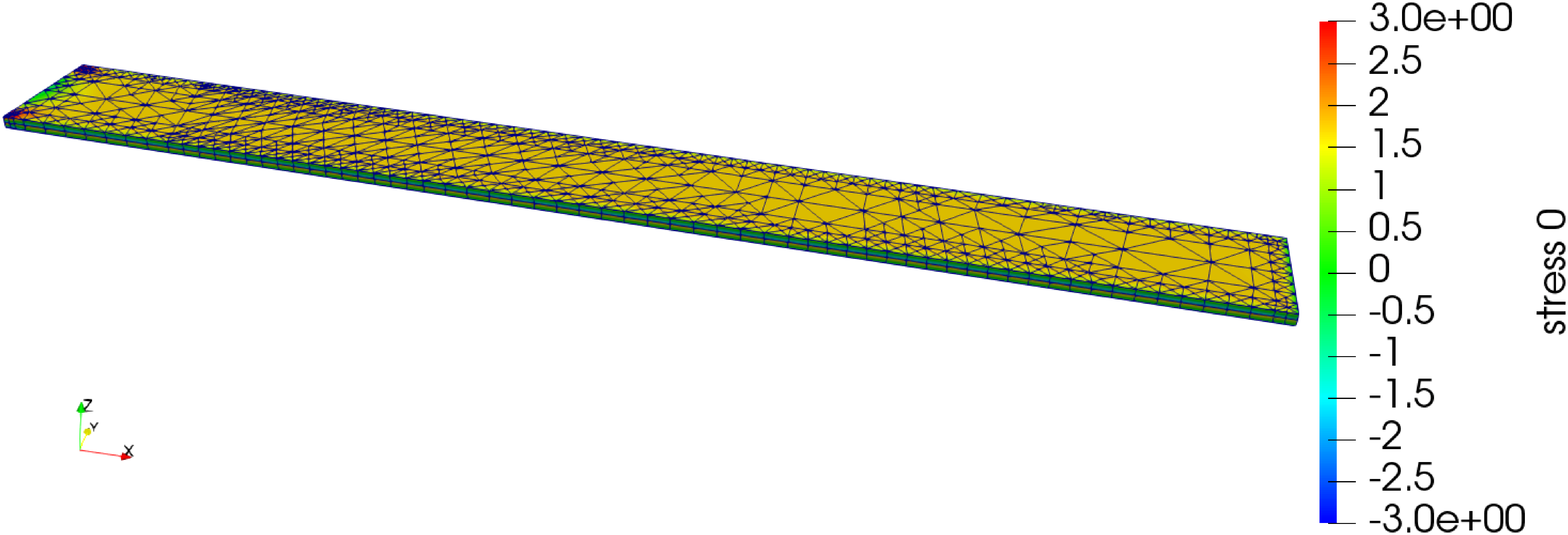}
\end{center}
\caption{Stress component $\sigma_{xx}$ for the coarsest and most refined mesh for V2, $k=1$, $k_\phi = 0$. The adaptive refinement is done towards all edges of the beam, where singularities in the solution occur.} \label{fig:beamsigma}
\end{figure}\begin{figure}
\begin{center}
\includegraphics[width=0.8\textwidth]{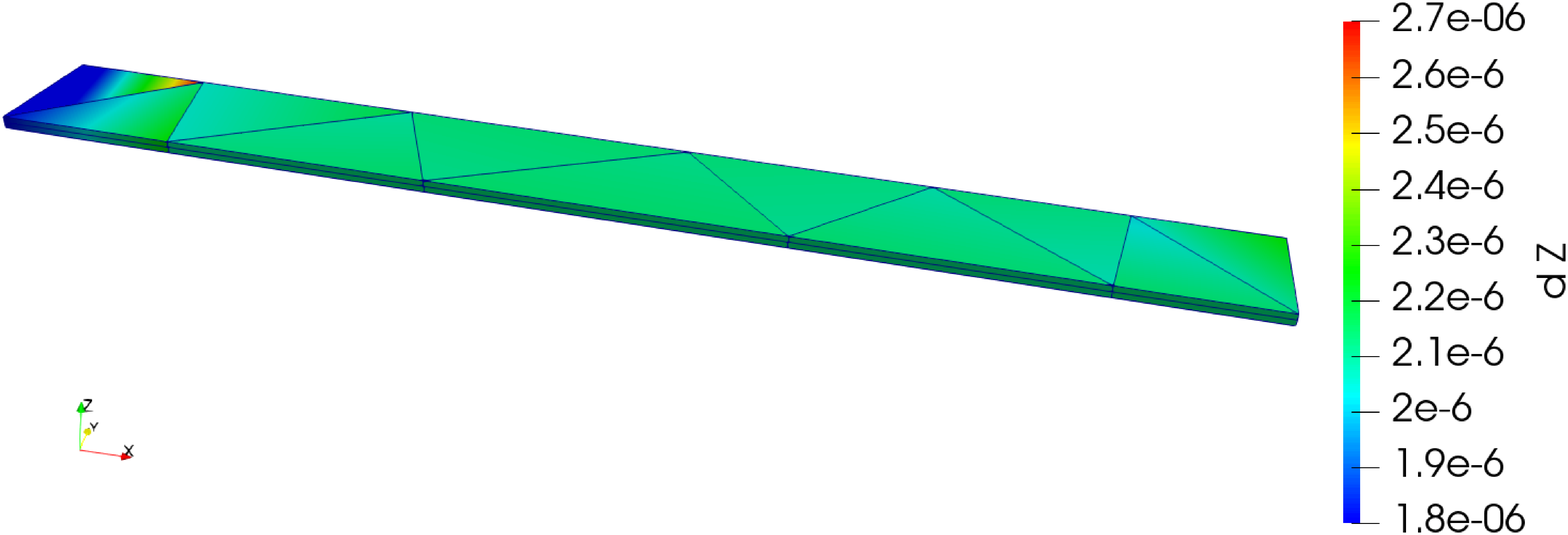}
\includegraphics[width=0.8\textwidth]{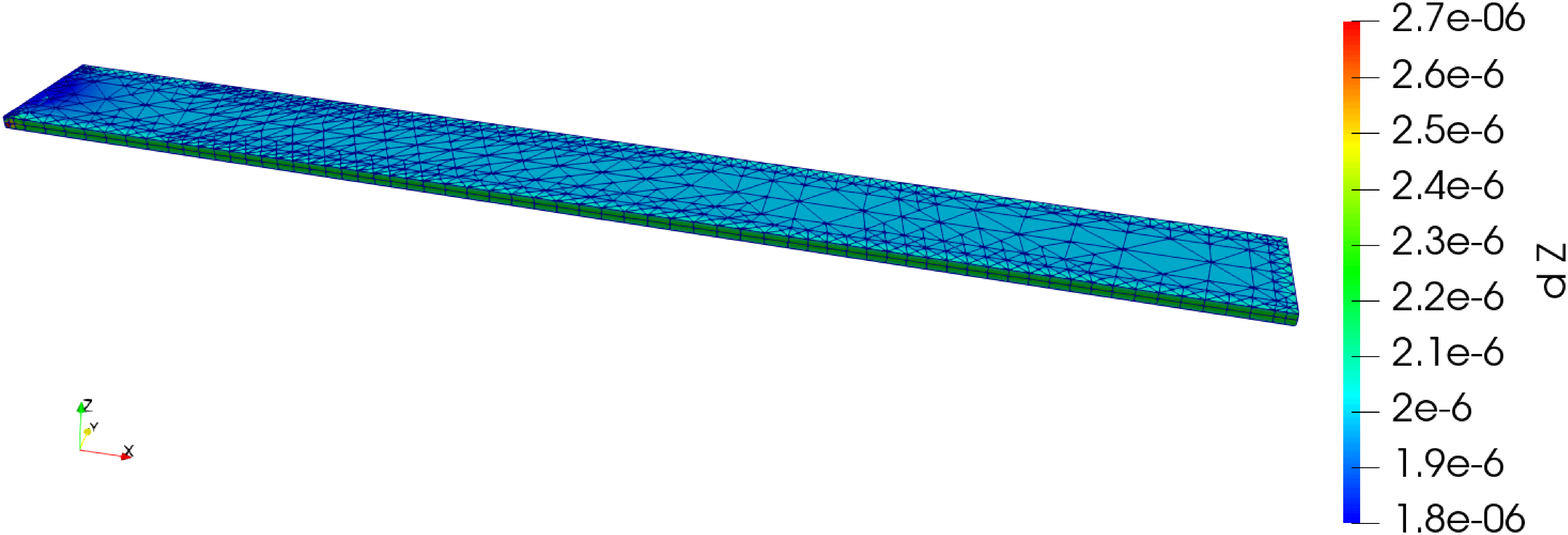}
\end{center}
\caption{Dielectric displacement component $D_z$ for the coarsest and most refined mesh for V2, $k=1$, $k_\phi = 0$. The adaptive refinement is done towards all edges of the beam, where singularities in the solution occur.} \label{fig:beamd}
\end{figure}

}

\begin{figure}
\begin{center}
\includegraphics[width=0.9\textwidth]{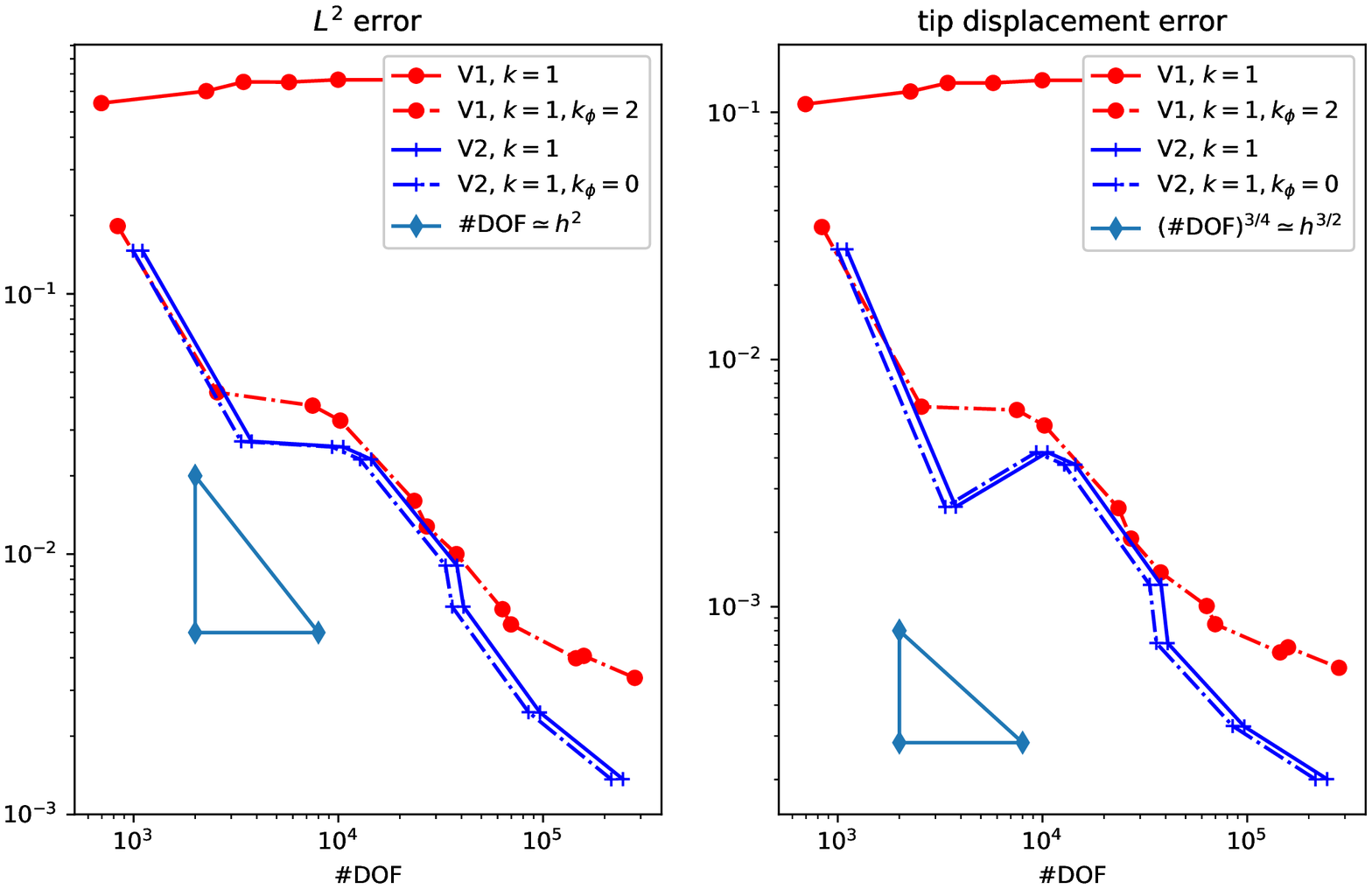}
\end{center}
\caption{Convergence of the displacments in $L^2$ sense $\|\uv - \uv_h\|_{L^2(\Omega)}$  and of average tip deflection  for  different methods with TDNNS order $k=1$ and potential order $k_\phi$. V1 indicates the first method using $\uv, \sigmat$ and $\phi$, while V2 indicates the second method using $\uv, \sigmat, \Dv$ and $\phi$.} \label{fig:bimbeam_conv1}
\end{figure}

\begin{figure}
\begin{center}
\includegraphics[width=0.9\textwidth]{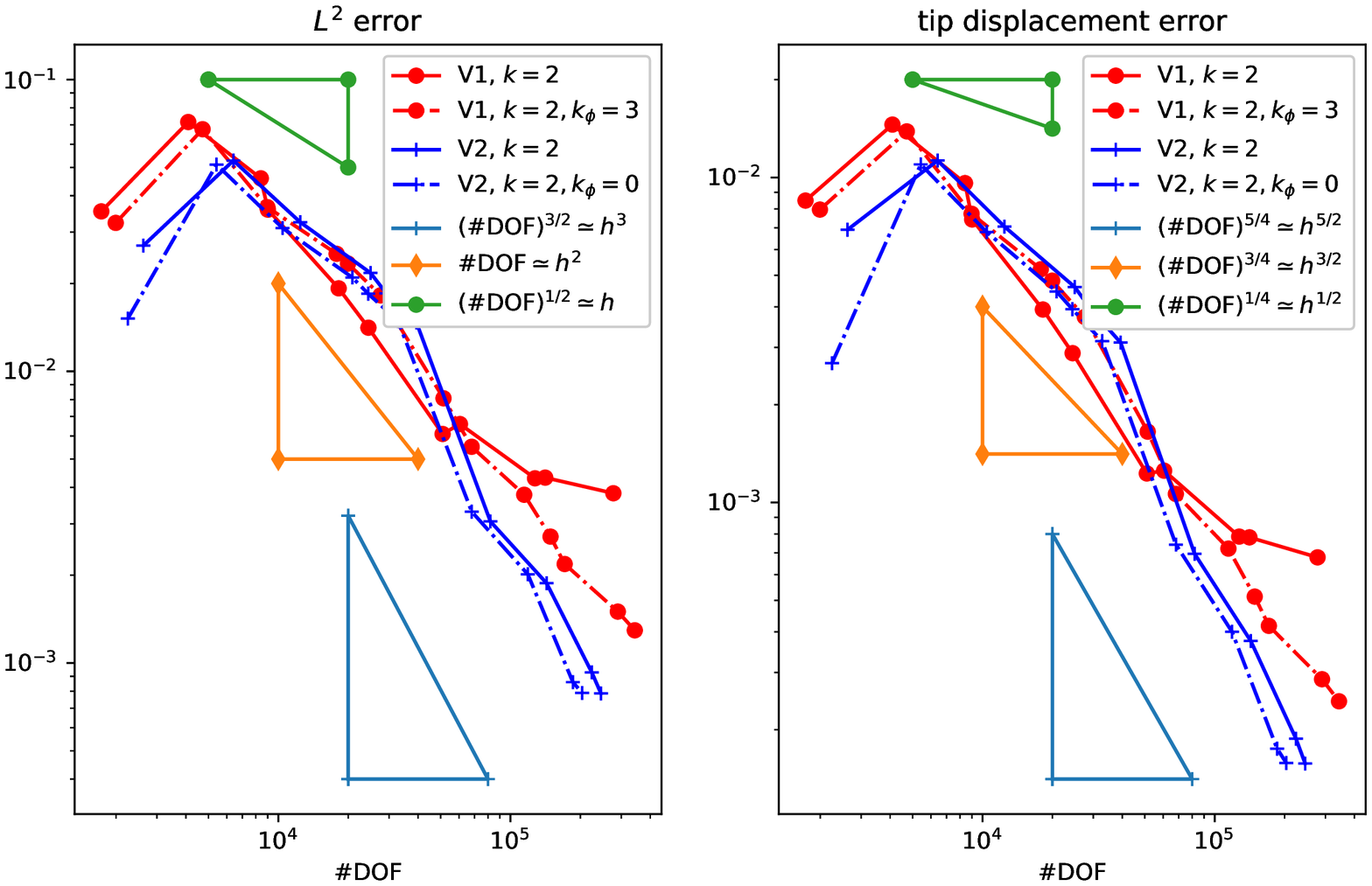}
\end{center}
\caption{Convergence of the displacments in $L^2$ sense $\|\uv - \uv_h\|_{L^2(\Omega)}$  and of average tip deflection for  different methods with  TDNNS order $k=2$ and potential order $k_\phi$. V1 indicates the first method using $\uv, \sigmat$ and $\phi$, while V2 indicates the second method using $\uv, \sigmat, \Dv$ and $\phi$.} \label{fig:bimbeam_conv2}
\end{figure}

\begin{figure}
\begin{center}
\includegraphics[width=0.9\textwidth]{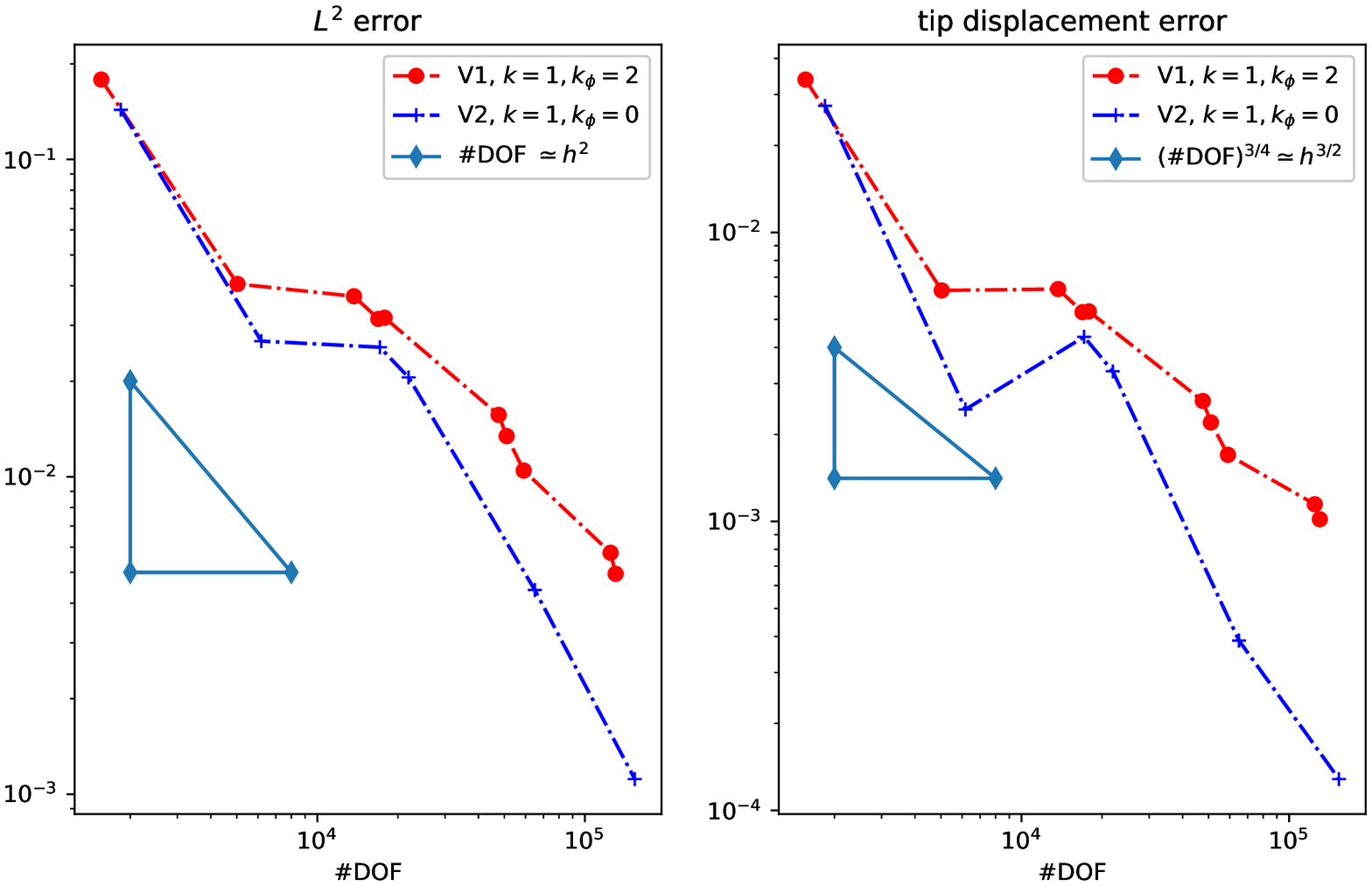}\\
\end{center}
\caption{Convergence of the displacments in $L^2$ sense $\|\uv - \uv_h\|_{L^2(\Omega)}$  and of average tip deflection  for  different methods with lowest TDNNS order $k=1$ and potential order $k_\phi$. Two elements per ply are used in thickness direction.} \label{fig:bimbeam_refz}
\end{figure}

In the first comparison, we use one element per ply in thickness direction. The results are depicted in Figure~\ref{fig:bimbeam_conv1} for TDNNS order $k=1$ and in Figure~\ref{fig:bimbeam_conv2} for TDNNS order $k=2$. For the lowest order approximations with $k=1$, we see that the convergence deteriorates, as soon as the error due to thickness discretization dominates. This deterioration is removed in a second comparison, where two elements per ply are used in thickness direction. The corresponding results are presented in Figure~\ref{fig:bimbeam_refz}.

We see that the lowest-order V1 approximation does not converge, since a linear electric potential element cannot recover the quadratic distribution of the potential. However, if we increase the order of the electric potential to $k_\phi = 2$ and leave the TDNNS order at $k=1$, we obtain convergence of order 2 in $L^2$ sense. The convergence rate deteriorates at approximately 200.000 unknowns, then the error due to the static thickness discretization dominates. The same optimal order of convergence is achieved by the second variant V2. Also here, the convergence rate deteriorates due to the static thickness discretization, but at a lower error level as for variant V1. The degrees of freedom are reduced, while the good approximation is preserved if we restrict the dielectric displacements to shape functions with constant divergence, see \eqref{eq:D0}--\eqref{eq:Phi0}. These results are indicated by $k_\phi = 0$. For the tip deflection, both V1 and V2 show convergence order $h^{3/2}$ initially.

For a higher-order TDNNS approximation $k=2$, we again get optimal convergence order 3 in $L^2$ sense for V2, and one order less for V1 with increased potential order $k_\phi=3$. Choosing the potential order $k_\phi = k = 2$, we see at best linear convergence in $L^2$ sense. For the tip deflection, all convergence orders are reduced by 1/2, as is to be expected for boundary evaluations.

In a second comparison, we use two elements per ply in thickness direction. We do computations for the lowest order elements of V1 and V2. We see that now the convergence order does not deteriorate, but is of optimal order 2 for the $L^2$ error, see Figure~\ref{fig:bimbeam_refz}.

At present, we do not aim at proofing any of these observed convergence orders mathematically. However, this shall be topic of further research.

\subsection{Homogenization of a $d_{15}$ MFC}

\begin{figure}
\begin{center}
\includegraphics[width=0.8\textwidth]{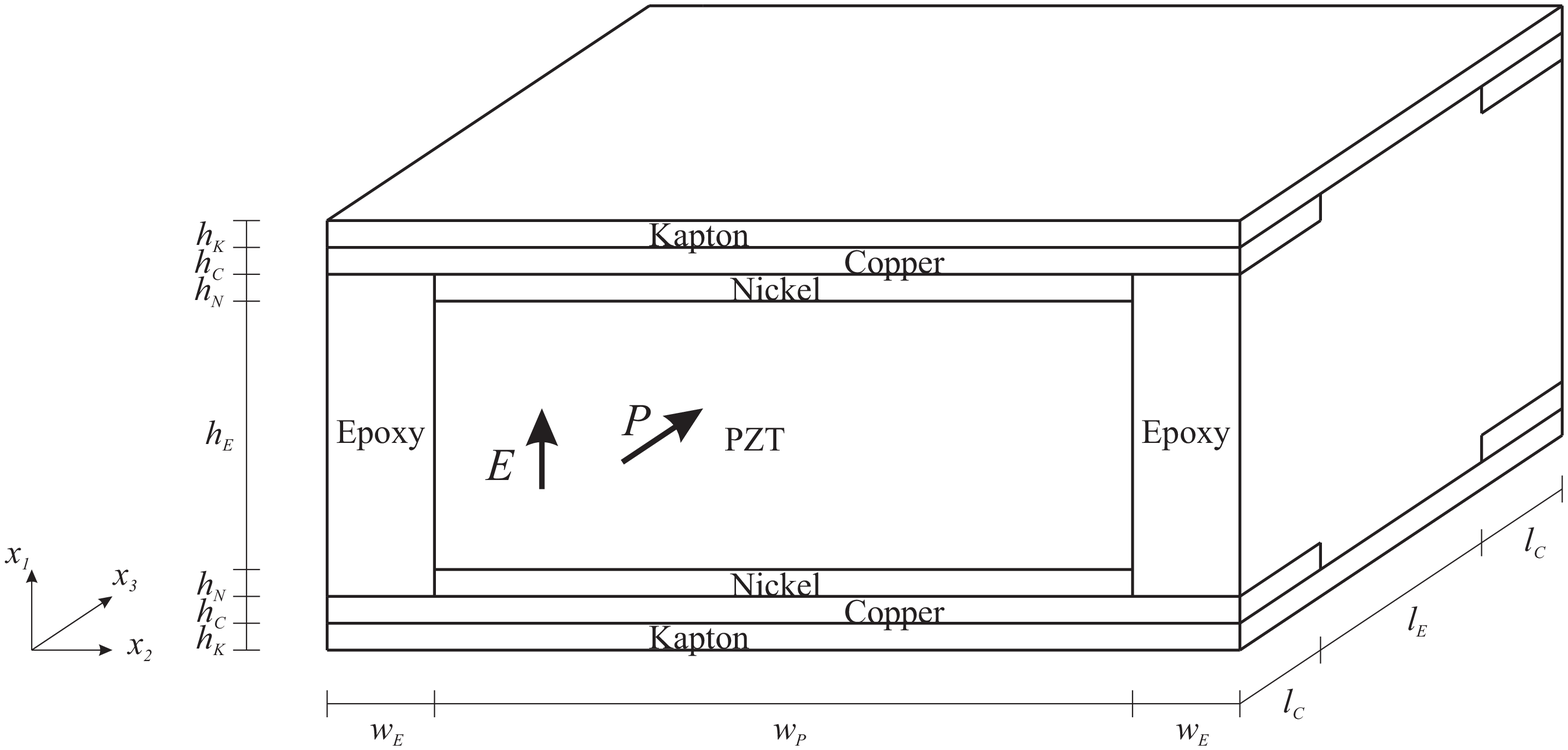}
\end{center}
\caption{Unit cell of the $d_{15}$ MFC. Dimensions: $h_P = \SI{180}{\micro\meter}$, $h_K = \SI{25}{\micro\meter}$, $h_N = \SI{2}{\micro\meter}$, $h_C = \SI{18}{\micro\meter}$, $l_E = \SI{420}{\micro\meter}$, $l_C = \SI{80}{\micro\meter}$, $w_E = \SI{27.5}{\micro\meter}$, $w_P = \SI{350}{\micro\meter}$.} \label{fig:MFC}
\end{figure}

As a second example, we compute some effective material properties of  a $d_{15}$ macro-fibre composite (MFC). We use the geometry and material data provided in \cite{KranzBenjeddouDrossel:2013,TrindadeBenjeddou:2013}. The setup is displayed in Figure~\ref{fig:MFC}. 

\begin{table}
\caption{Material constants used for the $d_{15}$ MFC.} \label{tab:MFCconstants}
\begin{center}
\begin{tabular}{|lr|lr|}
\hline
\multicolumn{4}{|l|}{SONOX P502}\\\hline
$S^{E}_{11}$ & $18.5\times 10^{-12}$ m$^2$/N\qquad \qquad &  $d_{31}$ & $-1.85\times10^{-10}$ m/V  \\
$S^{E}_{33}$ & $20.7\times 10^{-12}$ m$^2$/N\qquad \qquad &  $d_{33}$ & $4.40\times10^{-10}$ m/V  \\
$S^{E}_{12}$ & $-6.29\times 10^{-12}$ m$^2$/N\qquad \qquad &  $d^{lin}_{15}$ & $5.60\times10^{-10}$ m/V  \\
$S^{E}_{13}$ & $-6.23\times 10^{-12}$ m$^2$/N\qquad \qquad &  $\epsilon^{\sigma,lin}_{11}/\epsilon_0$ & 1950  \\
$S^{E}_{44}$ & $33.2\times 10^{-12}$ m$^2$/N\qquad \qquad &  $\epsilon^{\sigma}_{33}/\epsilon_0$ & 1850 \\
$S^{E}_{66}$ & $52.3\times 10^{-12}$ m$^2$/N\qquad \qquad &   &\\\hline
\multicolumn{4}{|l|}{Epoxy}\\\hline
$Y$ & $2.5\times 10^{9}$ N/m$^2$\qquad \qquad &   $\epsilon/\epsilon_0$ & 4.25  \\
$\nu$ & $0.42$\qquad \qquad &  &   \\\hline
\multicolumn{4}{|l|}{Kapton}\\\hline
$Y$ & $2.5\times 10^{9}$ N/m$^2$\qquad \qquad &   $\epsilon/\epsilon_0$ & 3.4  \\
$\nu$ & $0.34$\qquad \qquad &  &   \\\hline
\multicolumn{4}{|l|}{Copper}\\\hline
$Y$ & $110\times 10^{9}$ N/m$^2$\qquad \qquad &   $\epsilon/\epsilon_0$ & 2000  \\
$\nu$ & $0.34$\qquad \qquad &  &   \\\hline
\multicolumn{4}{|l|}{Nickel}\\\hline
$Y$ & $200\times 10^{9}$ N/m$^2$\qquad \qquad &   
$\nu$ & $0.31$  \\\hline
\end{tabular}
\end{center}
\end{table}

The material parameters  are displayed in Table~\ref{tab:MFCconstants}. Note that the permittivity of the nickel electrodes is not needed, as the electric potential is assumed to take a given input voltage there. Therefore, only the mechanic deformation is computed on the electrodes, while they are excluded in the electric equations. The nickel electrodes were not regarded in \cite{TrindadeBenjeddou:2013}, but in \cite{KranzBenjeddouDrossel:2013}. However, we use the same homogenization techniques as proposed in the former reference, and compare our results to theirs. We see that the influence of the electrodes is very small, as the results match well.

\begin{figure}
\begin{center}
\includegraphics[width=0.45\textwidth]{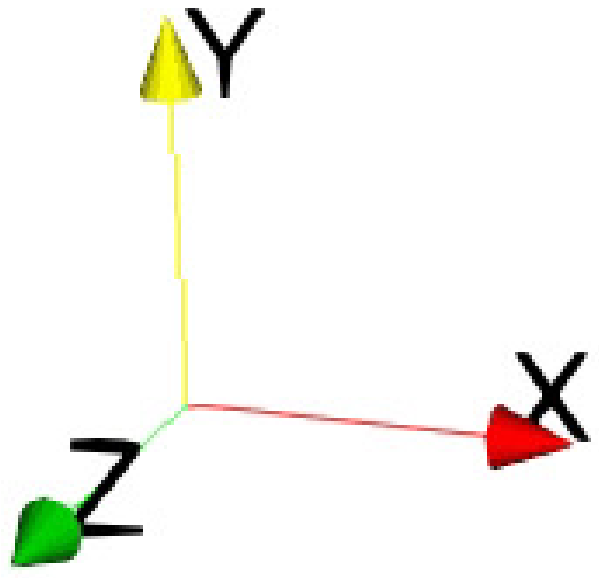}
\includegraphics[width=0.45\textwidth]{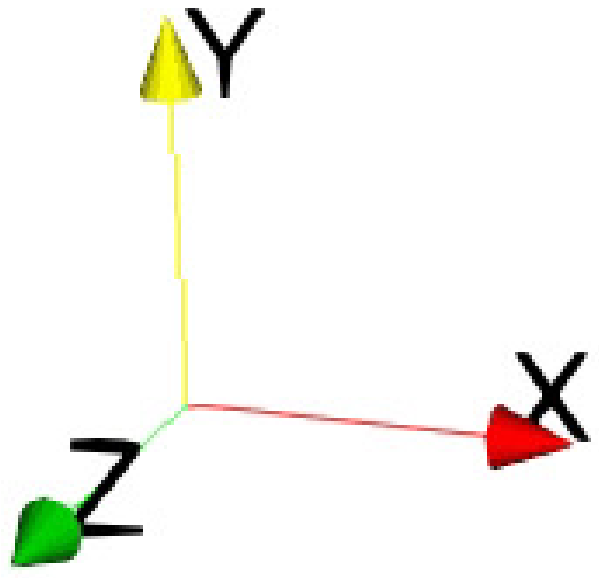}
\end{center}
\caption{Different finite element meshes used: coarse mesh with 560 elements (left) and fine mesh with 3164 elements (right).}\label{fig:meshes}
\end{figure}

We compute the shear modulus $G_{13} = 1/S^{E}_{55}$, the piezoelectric coefficient $d_{15}^{MFC}$ and dielectric constant $\epsilon_{11}^{\sigma,MFC}$ for the MFC using the proposed TDNNS finite element method. To this end, we implement local problems \#5 and \#7 from \cite{TrindadeBenjeddou:2013}. For the computation of $G^{MFC}_{13}$, periodic boundary conditions are prescribed \rev{for $\uv_t$ and $\sigma_nn$}, and an additional shear displacement is applied to the RVE. For the computation of $d_{15}^{MFC}$ and $\epsilon_{11}^{\sigma,MFC}$, the electric potential $\phi=1$V is prescribed on the electrodes, while stress-free conditions \rev{(i.e.\ $\sigma_{nn}=0$)} are assumed on the surface of the RVE. \rev{In variant V1, the electric potential is prescribed via its nodal values, in variant V2 it enters the right hand side of equation~\eqref{eq:potbc} on the surfaces of the electrode.}
From the finite element solutions of the respective load cases, the average shear strain and stress $\bar \gamma_{5}$ and $\bar \sigma_{5}$, the average electric field $\bar E_1$ and the average dielectric displacement $\bar D_1$ are computed by
\begin{align}
\bar \gamma_{5}& = \frac{2}{|\Omega|}\int_\Omega \eps_{13}\, d\Omega, &
\bar \sigma_{5}& = \frac{1}{|\Omega|}\int_\Omega \sigma_{13}\, d\Omega, &
\bar E_{1}& = \frac{1}{|\Omega|}\int_\Omega E_1\, d\Omega, &
\bar D_1& = \frac{1}{|\Omega|}\int_\Omega D_1\, d\Omega,.
\end{align}
Then the macroscopic piezoelectric and dielectric constants can be evaluated by
\begin{align} \label{eq:effectiveconstants}
G_{13}^{MFC} &= \bar \sigma_{5} / \bar \gamma_{5}, &
d_{15}^{MFC} &= \bar \gamma_{5} / \bar E_1, &
\epsilon_{11}^{\sigma,MFC} &= \bar D_1 / \bar E_1.
\end{align}
We use two different triangular finite element meshes in the $x_2x_3$ plane, which are extended to prismatic elements in $x_1$ direction. The electrode layer is always resolved by the finite element mesh. In thickness direction several different setups are considered: the coarsest using one element per ply (i.e. 7 layers in total), an intermediate setup using 4 elements in the active layer (i.e. 10 layers in total), and the finest using 8 elements in the active layer (i.e. 14 layers in total).
The in-plane and thickness direction are combined to three-dimensional tensor product meshes. The very coarsest and finest meshes are displayed in Figure~\ref{fig:meshes}.

\begin{figure}
\begin{center}
\includegraphics[width=0.46\textwidth]{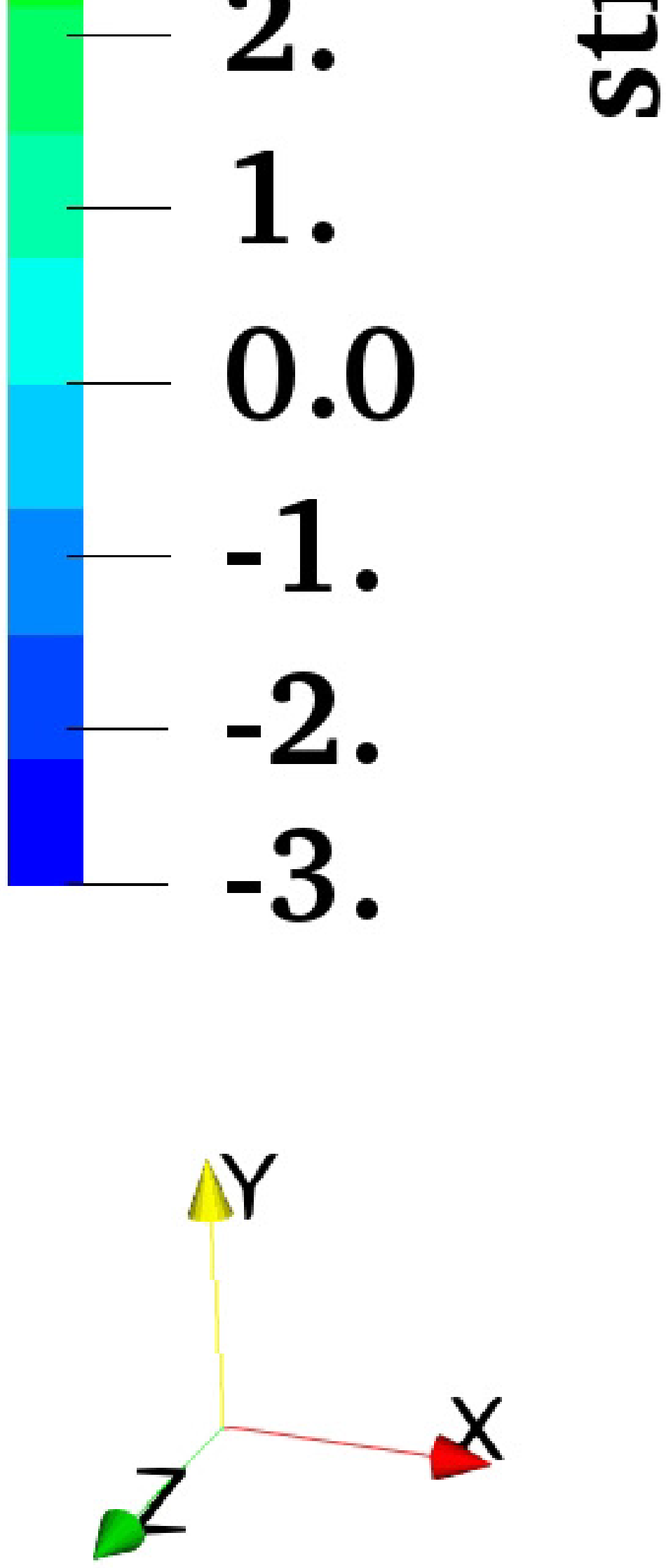}
\includegraphics[width=0.46\textwidth]{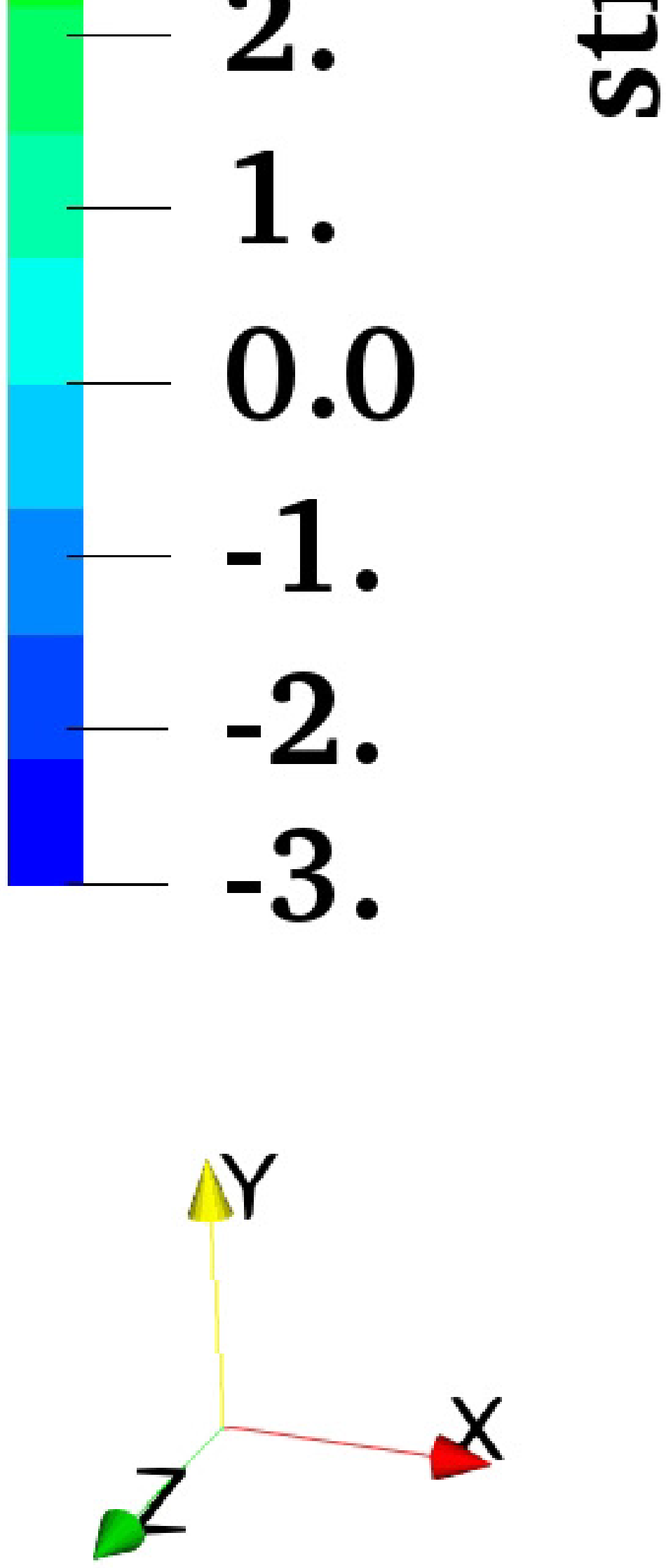}
\end{center}
\caption{Shear strain $\gamma_{5}$ from the homogenization problem for $G_{13}^{MFC}$, polynomial order $k=1$ on the coarsest mesh (left) and $k=2$ on the finest mesh (right), for TDNNS-based method V1.} \label{fig:strain5G}
\end{figure}

\begin{figure}
\begin{center}
\includegraphics[width=0.46\textwidth]{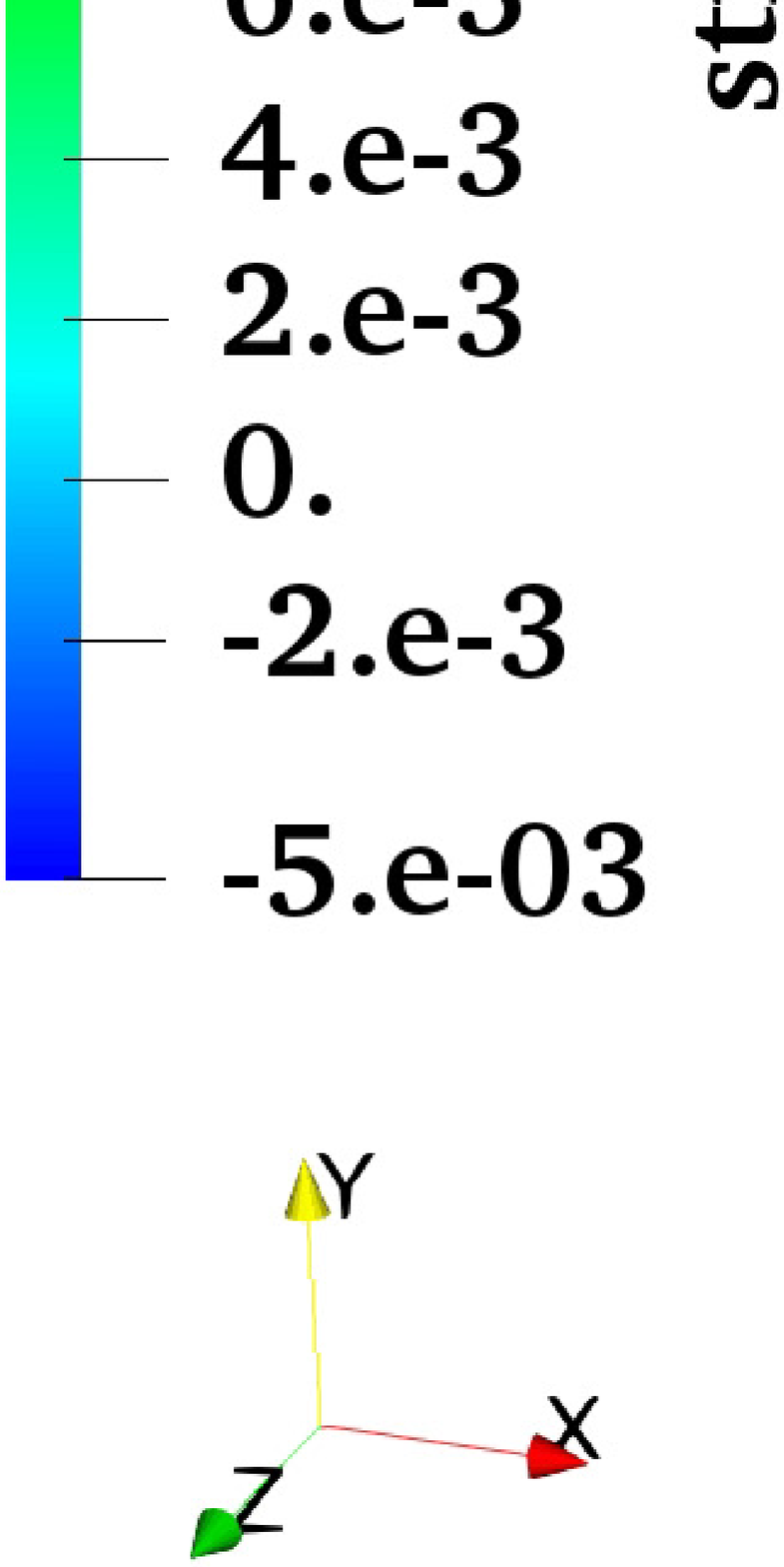}
\includegraphics[width=0.46\textwidth]{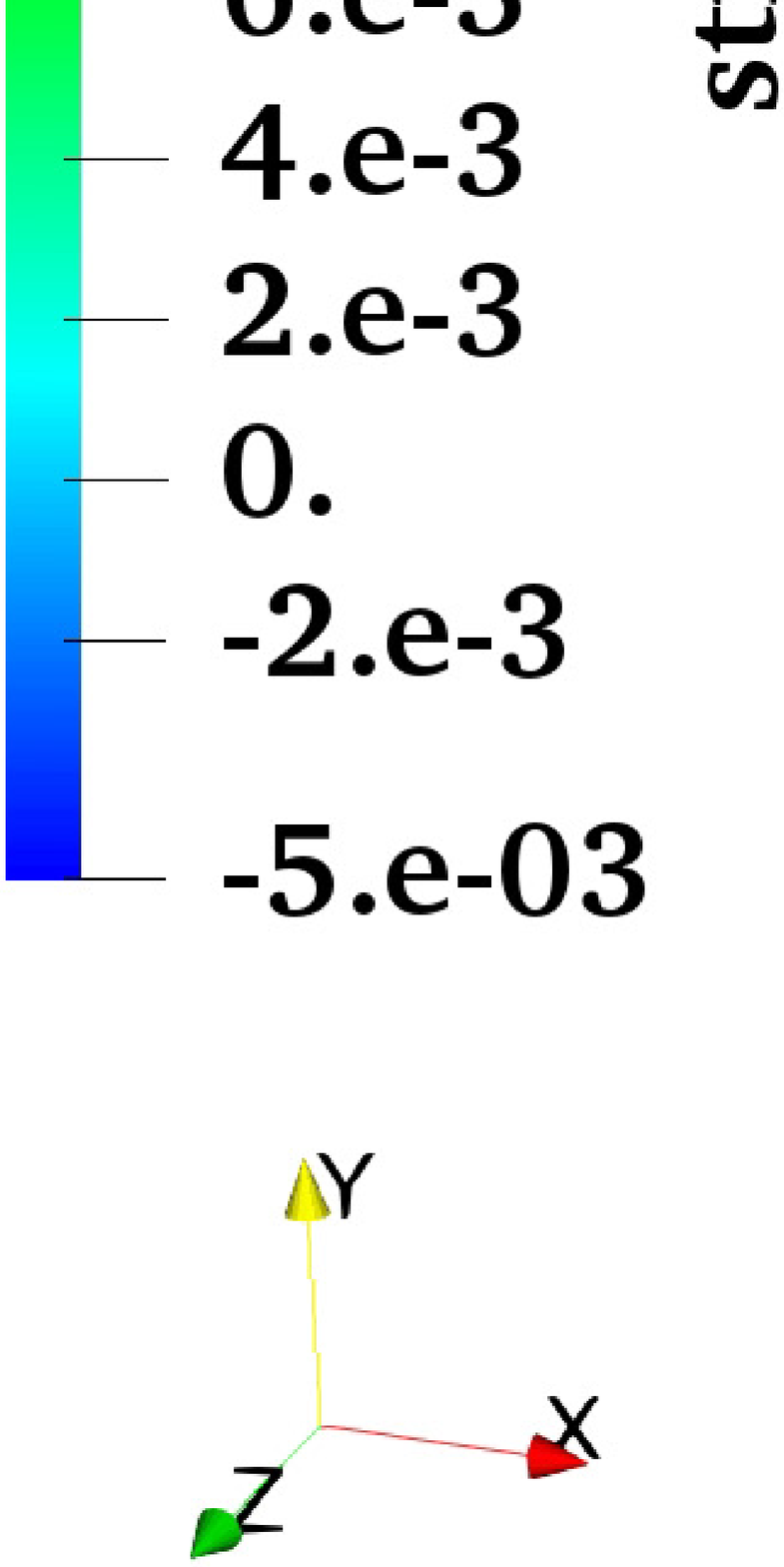}
\end{center}
\caption{Shear stress $\sigma_{5}$ from the homogenization problem for $G_{13}^{MFC}$, polynomial order $k=1$ on the coarsest mesh (left) and $k=2$ on the finest mesh (right), for TDNNS-based method V1.} \label{fig:stress5G}
\end{figure}

\begin{figure}
\begin{center}
\includegraphics[width=0.46\textwidth]{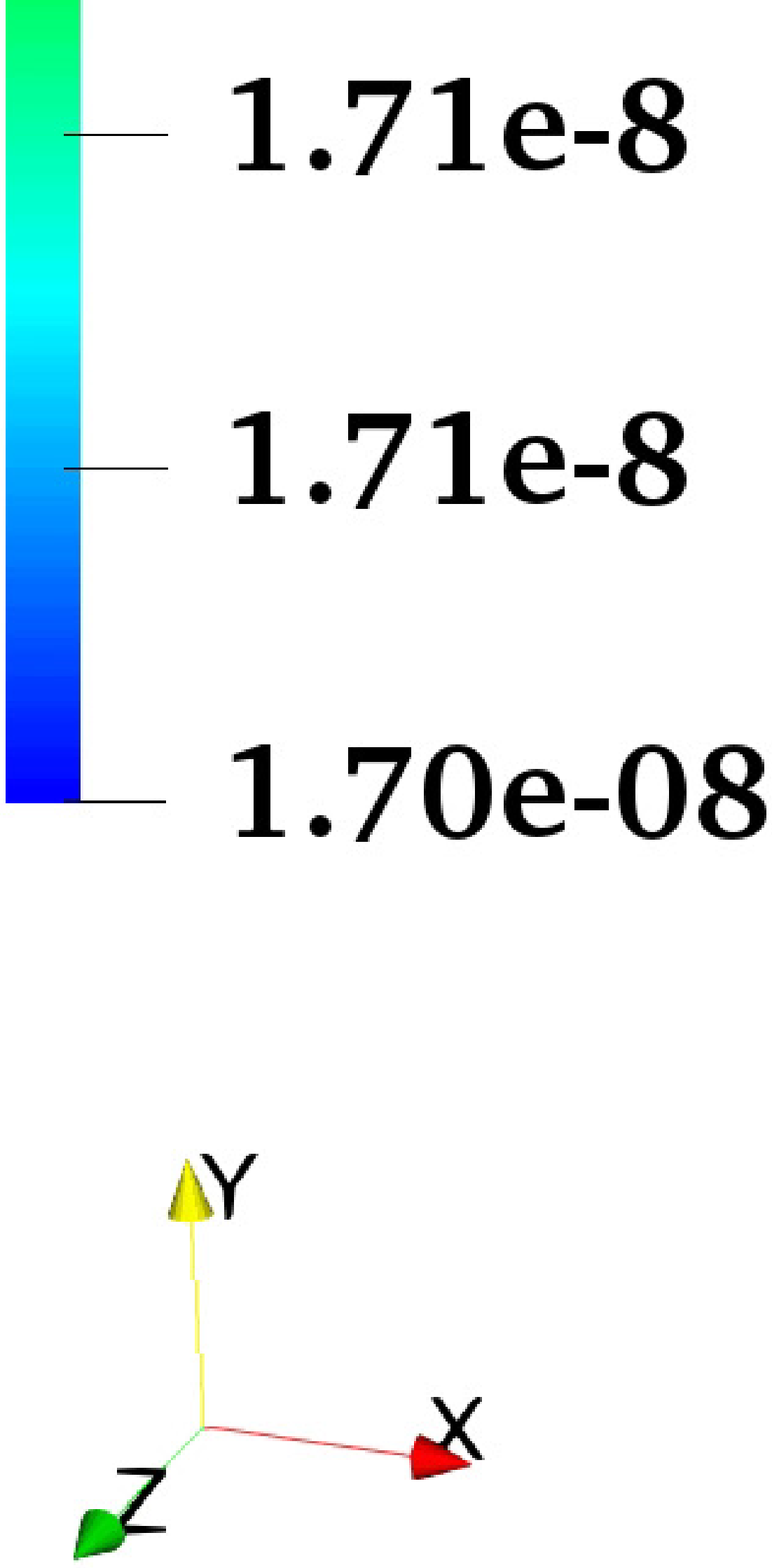}
\includegraphics[width=0.46\textwidth]{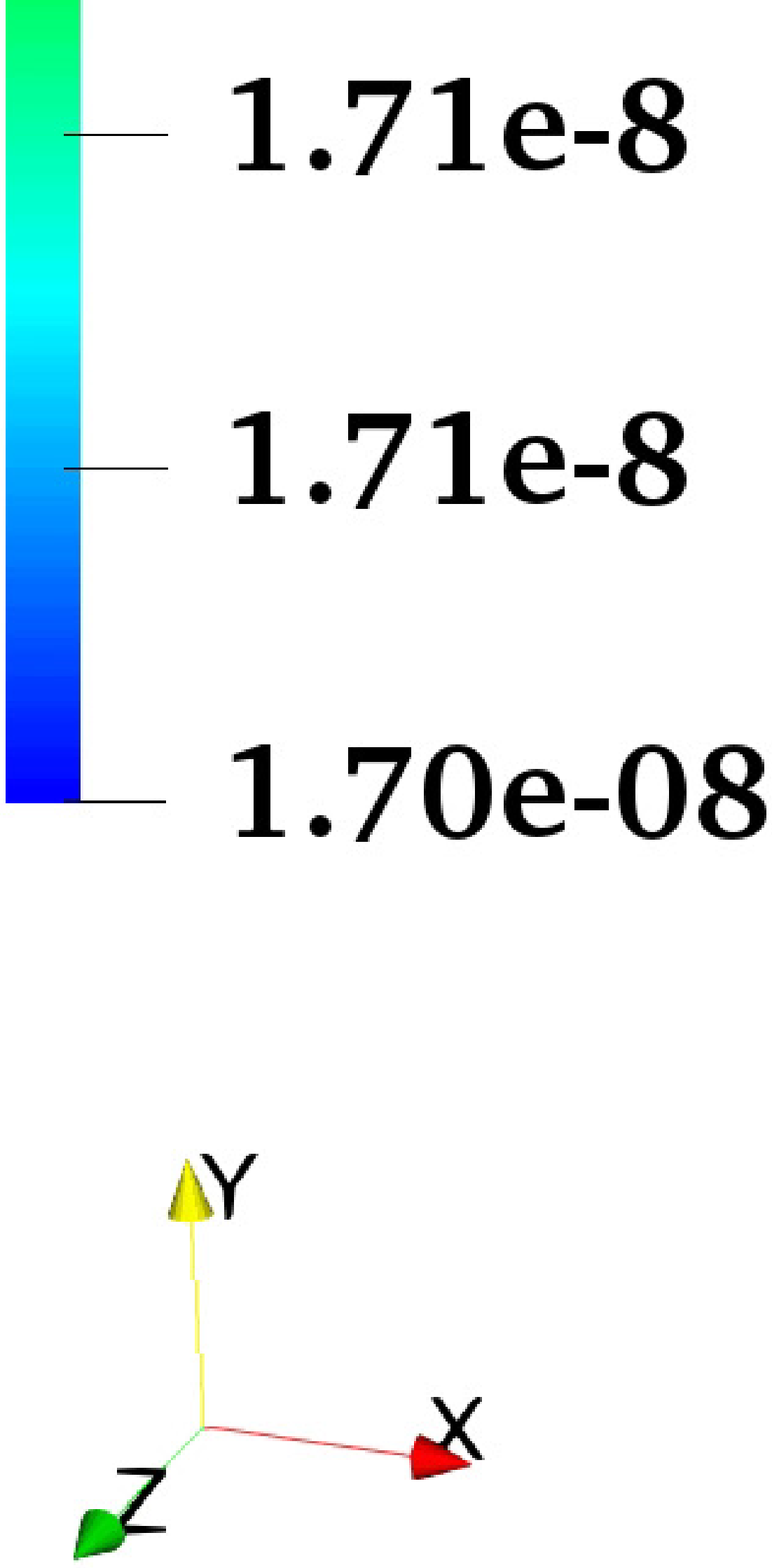}
\end{center}
\caption{Dielectric displacement $D_1$ from the homogenization problem for $d_{15}^{MFC}$ and $\epsilon_{11}^{\sigma,MFC}$, polynomial order $k=1$ on the coarsest mesh (left) and $k=3$ on the finest mesh (right), for TDNNS-based method V1.}
\label{fig:D1d}
\end{figure}

\begin{figure}
\begin{center}
\includegraphics[width=0.46\textwidth]{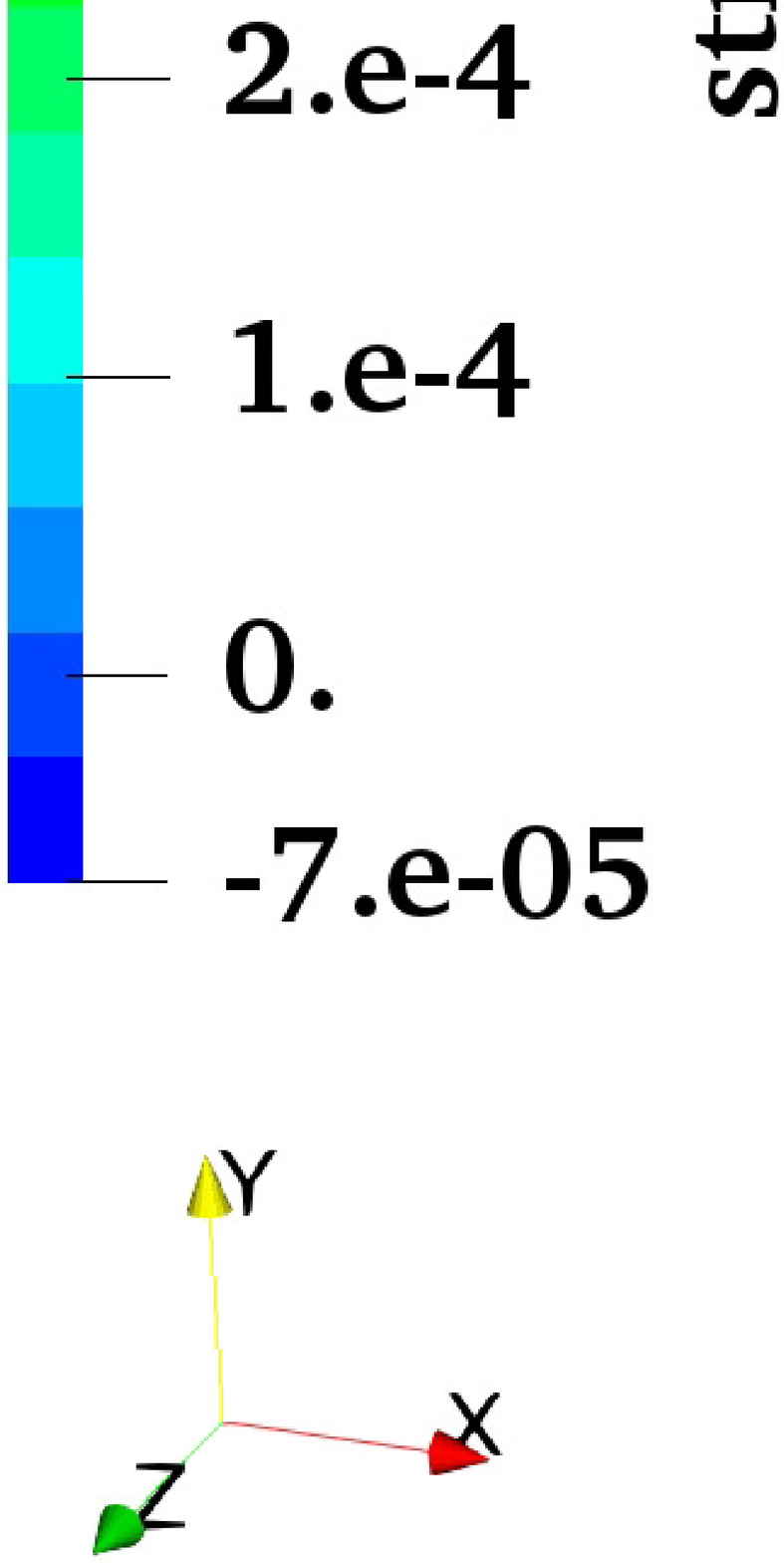}
\includegraphics[width=0.46\textwidth]{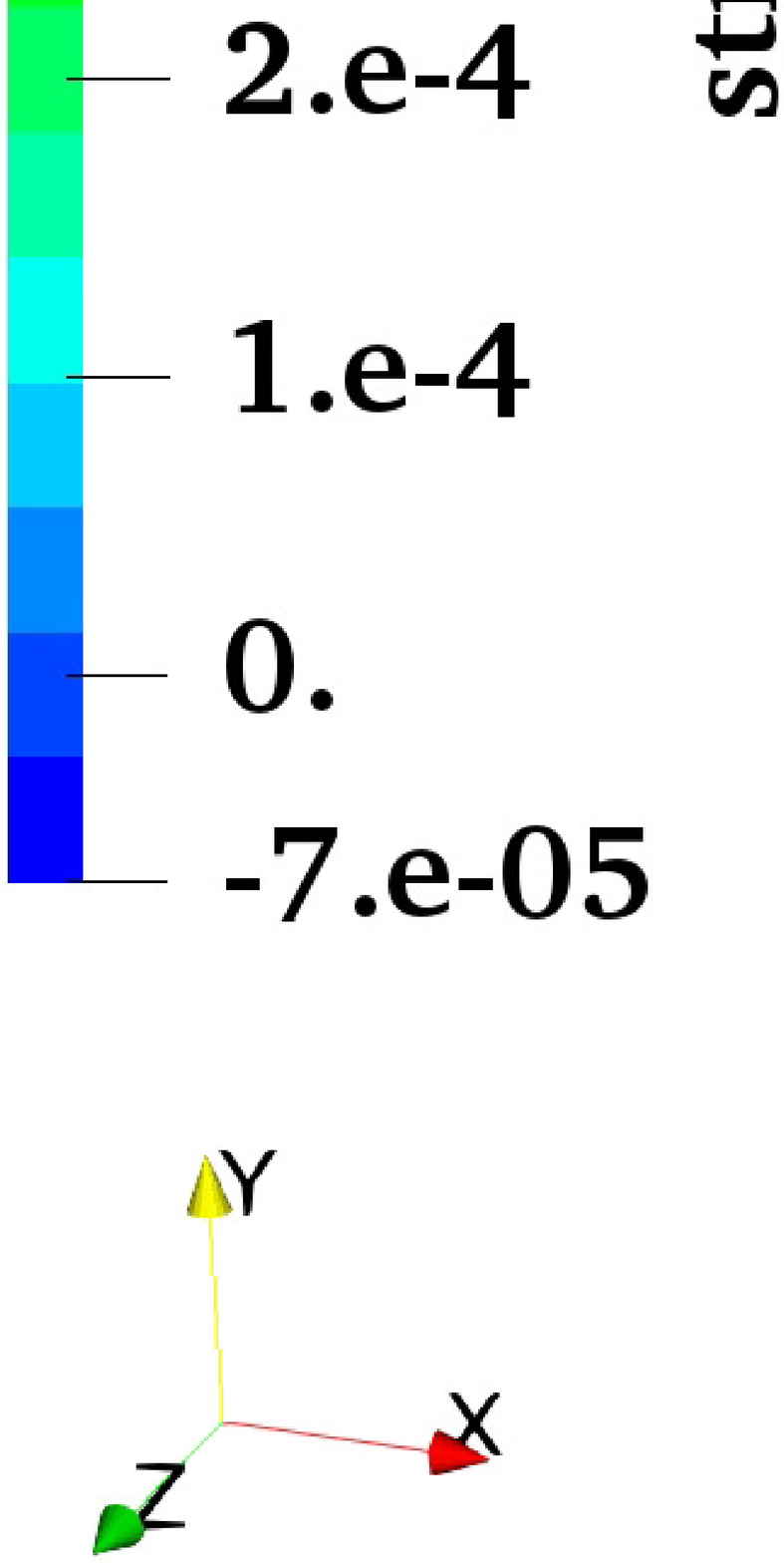}\\
\end{center}
\caption{Shear strain $\gamma_{5}$ from the homogenization problem for $d_{15}^{MFC}$ and $\epsilon_{11}^{\sigma,MFC}$, polynomial order $k=1$ on the coarsest mesh (left) and $k=3$ on the finest mesh (right), for TDNNS-based method V1.} \label{fig:strain5d}
\end{figure}

Effective values are computed from formulae \eqref{eq:effectiveconstants} for the different finite element meshes, variants V1 and V2 and different polynomial orders. The results are collected in Table~\ref{tab:linearcoefG13}, Table~\ref{tab:linearcoefd15} and Table~\ref{tab:linearcoefeps11}. They are compared to values provided in \cite{TrindadeBenjeddou:2013}. One can see that the values are close, the small differences may arise from the fact that we resolve the nickel electrodes by the finite element mesh. Figure~\ref{fig:strain5G} to Figure~\ref{fig:strain5d} show various computed fields for the electric potential based method V1. In all figures, the left hand plot shows the field for the coarsest mesh at lowest polynomial order, while the right hand plot shows the solutions on the finest mesh at higher polynomial order.

\begin{table} 
 \caption{Effective shear modulus $G^{MFC}_{13}$ of the MFC, computed by a homogenization method, compared to the values provided in \cite{TrindadeBenjeddou:2013} . } \label{tab:linearcoefG13}
\begin{center}
\begin{tabular}{|l|c|c|c|c|c|c|}\hline
 & \multicolumn{3}{l|}{V1} & \multicolumn{3}{l|}{V2}  \\
$G_{13}^{MFC}$, [GPa]& $k=1$ & $k=2$ & $k=3$ & $k=1$ & $k=2$ & $k=3$\\ \hline 
coarse mesh & & & & & & \\
7 layers & 3.141 & 3.131 & 3.139 & 3.142 & 3.132 & 3.140\\
10 layers & 3.135 & 3.136 & 3.137 & 3.135 & 3.136 & 3.138\\
14 layers & 3.133 & 3.136 & 3.138 & 3.144 & 3.136  & 3.138 \\
fine mesh & & & & & & \\
7 layers & 3.154 & 3.135 & 3.138 & 3.156 & 3.136& 3.138\\
10 layers & 3.144 & 3.137 & 3.137 & 3.144 & 3.137 & 3.137\\
14 layers & 3.144 & 3.137 & 3.138 & 3.134 & 3.137  & 3.138 \\
\hline\hline
\multicolumn{5}{|l|}{Trindade and Benjeddou 2013} &  \multicolumn{2}{l|}{3.10}\\
\hline
\end{tabular}
\end{center}
\end{table}

\begin{table} 
 \caption{Effective piezoelectric coefficient $d^{MFC}_{15}$ of the MFC, computed by a homogenization method, compared to the values provided in \cite{TrindadeBenjeddou:2013}. } \label{tab:linearcoefd15}
\begin{center}
\begin{tabular}{|l|c|c|c|c|c|c|}\hline
& \multicolumn{3}{l|}{V1} & \multicolumn{3}{l|}{V2}  \\
$d_{15}^{MFC}$, [pC/N] & $k=1$ & $k=2$ & $k=3$ & $k=1$ & $k=2$ & $k=3$\\ \hline 
coarse mesh & & & & & & \\
7 layers &  554.001 & 553.470 & 553.515 & 556.992 & 554.646 & 553.869\\
10 layers & 554.428 & 553.637 & 553.549 & 555.274 & 553.948 & 553.708\\
14 layers & 554.477 & 553.636 & 553.548 & 555.267 & 553.935  &553.702\\
\hline
fine mesh & & & & & & \\
7 layers  & 553.972 & 553.424 & 553.513 & 556.644 & 554.462 & 553.929\\
10 layers & 554.345 & 553.618 & 553.541 & 554.839 & 553.778 & 553.616\\
14 layers & 554.357 & 553.620 & 553.542 & 554.791 & 553.773  & 553.616 \\
\hline\hline
\multicolumn{5}{|l|}{Trindade and Benjeddou 2013} &  \multicolumn{2}{l|}{554.02}\\
\hline
\end{tabular}
\end{center}
\end{table}

\begin{table} 
 \caption{Effective dielectric constant $\epsilon^{\sigma,MFC}_{11}$, computed by a homogenization method, compared to the values provided in \cite{TrindadeBenjeddou:2013}} \label{tab:linearcoefeps11}
\begin{center}
\begin{tabular}{|l|c|c|c|c|c|c|}\hline
& \multicolumn{3}{l|}{V1} & \multicolumn{3}{l|}{V2}  \\
$\epsilon_{11}^{\sigma,MFC}$, [nF/m] & $k=1$ & $k=2$ & $k=3$ & $k=1$ & $k=2$ & $k=3$\\ \hline 
coarse mesh & & & & & & \\
7 layers & 14.962 & 15.001 & 14.950 & 15.041 & 15.0315 & 15.012\\
10 layers & 14.999 & 15.005 & 14.955 & 15.022 & 15.013 & 15.010\\
14 layers & 14.999 & 15.005 & 14.955 & 15.020 & 15.013  & 15.010 \\
fine mesh & & & & & & \\
7 layers & 14.963 & 15.001 & 15.002 & 15.033 & 15.028 & 15.011\\
10 layers & 15.001 & 15.006 & 15.007 & 15.015 & 15.010 & 15.009\\
14 layers & 15.002 & 15.006 & 15.007 & 15.013 & 15.010  & 15.009 \\
\hline\hline
\multicolumn{5}{|l|}{Trindade and Benjeddou 2013} &  \multicolumn{2}{l|}{15.02}\\
 \hline
\end{tabular}
\end{center}
\end{table}

\section{Appendix -- Implementation in Netgen/NGSolve}

As already mentioned, the finite elements above are implemented in the open-source finite element software package Netgen/NGSolve. In the sequel, we show how to set up such a finite element problem in a python script.
To this end, we use the example of a piezoelectric bimorph beam from the previous section.

\begin{figure}
\begin{center}
\includegraphics[width=0.6\textwidth]{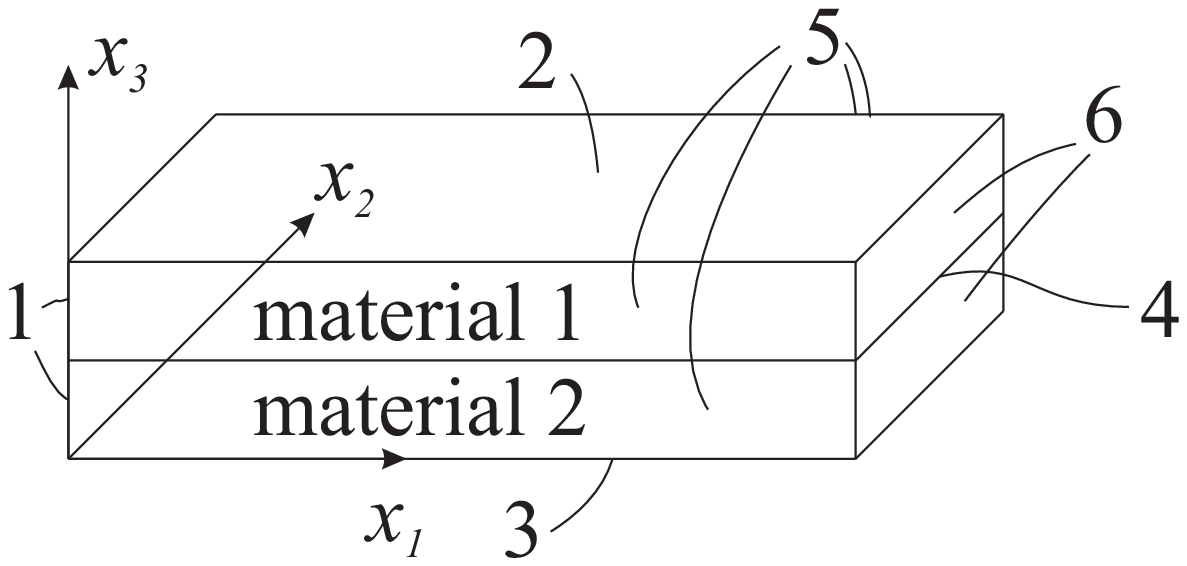}
\end{center}
\caption{Numbering of material and boundary regions of the bimorph beam.} \label{fig:bimbeam_regions}
\end{figure}

We present the essential steps of the implementation of the electric potential based method and the dielectric displacement based method. First, we define the geometry and generate the prismatic finite element mesh. The different material and boundary regions are numbered as displayed in Figure~\ref{fig:bimbeam_regions}. For more information on geometry and mesh generation in two and three space dimensions see the chapter on geometric modelling and mesh generation at \verb+https://ngsolve.org/docu/nightly/i-tutorials/+.
To define the geometry of the bimorph beam, constructive solid geometry CSG is used. The beam is defined as the intersection of half spaces. Each half space is given by \verb+Plane(p, n)+, where the tuple \verb+p+ is a point on the plane and \verb+n+ is the outer normal vector. The boundary condition number is set by \verb+bc+.
\begin{verbatim}
from netgen.csg import *
from ngsolve import *

p_left =  Plane(Pnt(0,0,0),  Vec(-1,0,0)).bc(1)
p_right = Plane(Pnt(l,0,0),  Vec(1,0,0)) .bc(6)
p_top =   Plane(Pnt(0,0,h),  Vec(0,0,1)) .bc(2)
p_bottom =  Plane(Pnt(0,0,-h), Vec(0,0,-1)).bc(3)
p_center = Plane(Pnt(0,0,0),  Vec(0,0,1)) .bc(4)
p_front =  Plane(Pnt(0,b,0),  Vec(0,1,0)) .bc(5)
p_back = Plane(Pnt(0,0,0),  Vec(0,-1,0)).bc(5)

geometry = CSGeometry()
matnr_1 = geometry.Add((p_left * p_right * p_top * p_front * p_back) - p_center)
matnr_2 = geometry.Add(p_left * p_right* p_bottom * p_center * p_front * p_back)
\end{verbatim}

To generate a prismatic mesh, the electroded surfaces have to be identified. If more than one element per ply in thickness direction should be used, positions of slices can be given and a \verb+ZRefinement+ called.

\begin{verbatim}
geometry.CloseSurfaces(p_center, p_bottom, slices=[0.5])
geometry.CloseSurfaces(p_center, p_top, slices = [0.5])

netgenmesh = geometry.GenerateMesh(maxh=100)
ZRefinement(netgenmesh, geometry)  # optional
mesh = Mesh(netgenmesh)
\end{verbatim}

\subsection{Electric potential based method}

First, we define coefficient functions that resemble the material parameters $\St^{E}$, $\dt$ and $\epsilont^\sigma$. Assuming we have tuples of length 36, 18 and 9 containing the respective constant material parameters, matrix-valued coefficient functions are defined by
\begin{verbatim}
# SE_tup = (S11, S12, ... S66)
# d_tup = (d11, d12, ... d36)
# epsilonT_tup = (epsT11, epsT12, ... epsT33)
SE = CoefficientFunction(SE_tup , dims=(6,6))
d = CoefficientFunction(d_tup, dims=(3,6))
epsilonT = CoefficientFunction(epsilonT_tup , dims=(3,3))
\end{verbatim}
The unit outward normal on boundaries and element interfaces is often needed in the TDNNS method. It is available in NGSolve as a special coefficient function
\begin{verbatim}
n = specialcf.normal(3)
\end{verbatim}

For the electric potential based method, we need three different finite element spaces: the stress space $\spaceSigma_h$, the displacement space $\spaceV_h$ and the potential space $\spacePhi_h$. We collect these three spaces into one compound space $\spaceX_h$, where they are ordered consecutively, by
%
\begin{verbatim}
Sigma = HDivDiv(mesh, order=k, dirichlet=[2,3,5,6] )
V = HCurl(mesh, dirichlet=[1], order=k)
Phi = H1(mesh, order=k_phi, dirichlet=[2,3,4] )
X = FESpace([Sigma, V, Phi])
\end{verbatim}
%
The keyword \verb+dirichlet+ marks boundary regions, where essential boundary conditions on the normal stress, tangential displacements and electric potential are enforced.

Next, the global solution vector \verb+U+ containing $\sigmat, \uv$ and $\phi$ is defined. In the current example, we have an inhomogeneous boundary condition for the electric potential at the outer electrodes. To implement this boundary condtion, we split the electric potential in two parts,
\begin{equation}
\phi = \tilde \phi + \phi_0.
\end{equation}
The second part $\phi_0$ satisfies the non-zero boundary condition $\phi_0 = 75 V$ at the electrodes, and is set in advance. The first part $\tilde \phi$ satisfies the homogeneous boundary conditions $\tilde \phi = 0$ at all electrodes, and is computed by the finite element method. This is realized in the python code as
\begin{verbatim}
U = GridFunction(X)
U0 = GridFunction(X)
Stress, Disp, Pot = U.components
Pot_0 = U0.components[2]
Pot_0.Set([0,75,75,0,0,0,0], VOL_or_BND=BND)
\end{verbatim}

Finally, the variational equations have to be defined in symbolic form. To this end, we introduce (symbolic) trial and test functions resembling $\sigmat, \uv, \tilde\phi$ and $\delta \sigmat, \delta \uv, \delta \phi$.
\begin{verbatim}
sigma, u, tilde_phi = X.TrialFunction()
d_sigma, d_u, d_phi = X.TestFunction()
\end{verbatim}
A few definitions that are useful to make the code more readable are given below. While the first function computing the tangential component of a vector is obvious, we mention that the second function gives the stress tensor in six-dimensional engineering vector notation. The last function is the divergence, which is the pre-implemented derivative for normal-normal continuous \verb+HDivDiv+ functions.
\begin{verbatim}
def tang(u): return u - InnerProduct(u,n)*n
def vec(sigma): return sigma.Operator("vec")
def div(sigma): return sigma.Deriv()
\end{verbatim}
The left hand side of eq.~\eqref{eq:V1} is summarized in the bilinear form \verb+a+, while the right hand side is represented by the linear form \verb+f+. The bilinear form produces the stiffness matrix, the linear form the load vector.
\begin{verbatim}
a = BilinearForm(X)
a += SymbolicBFI( InnerProduct (vec(sigma), SE*vec(d_sigma) ) )
a += SymbolicBFI(-InnerProduct (d*vec(sigma), d_phi.Deriv()) \
                 -InnerProduct(d*vec(d_sigma), tilde_phi.Deriv()) )
a += SymbolicBFI( InnerProduct(tilde_phi.Deriv(),epsilonT*d_phi.Deriv()))
a += SymbolicBFI( InnerProduct(div(sigma),d_u)+InnerProduct(div(d_sigma), u))
a += SymbolicBFI(-InnerProduct(sigma*n,tang(d_u))\
                 -InnerProduct(d_sigma*n,tang(u)), element_boundary=True)

f = LinearForm(X)
f += SymbolicLFI(-InnerProduct(Pot_0.Deriv(), epsilonT*d_phi.Deriv()) )
f += SymbolicLFI( InnerProduct(d*vec(d_sigma), Pot_0.Deriv()) )
\end{verbatim}

After these definitions, stiffness matrix and load vector are assembled, the inverse of the stiffness matrix is computed and applied to the load vector. The complete solution vector is computed by adding $\phi_0$. Paraview output is generated, and the tip displacement is evaluated in two different ways.
\begin{verbatim}
a.Assemble()
f.Assemble()

invmat = a.mat.Inverse(X.FreeDofs(), inverse="umfpack")     
U.vec.data = invmat * f.vec
U.vec.data += U0.vec

vtk = VTKOutput(ma=mesh,coefs=[vec(Stress), Disp, Pot, -Pot.Deriv()],
names=["stress","disp","Phi", "E"],
filename='BimorphBeam', subdivision=2)
vtk.Do()

bar_uz = Integrate(Disp[2], mesh, BND, region_wise=True)
print("av. tip displacement = ", bar_uz[5]/(2*h*b))

uz = Disp(mesh(l,b/2,0))[2]
print("displacement at point (l, b/2, 0) = ", uz)
\end{verbatim}

To avoid numerical problems, the solution of the linear system should be done in several steps. The interior degrees of freedom are eliminated from the global system by static condensation, and computed by a local postprocessing steps. This leads to better conditioned systems, for more details on the implementation see \verb+https://ngsolve.org/docu/latest/how_to/howto_staticcondensation.html+.

\subsection{Dielectric displacement based method}

For the implementation of the dielectric displacement based method, many steps are identical or similar to the potential based method. Thus we concentrate on those steps that differ.

We assume the material constants are given as coefficient functions in $\gt$-form. The electric potential is fully discontinuous, marked by the keyword \verb+L2+. The finite elements consist of four parts now, as the dielectric displacements are added as further unknown. The dielectric displacements are modeled in the normal-continuous space \verb+HDiv+.  However, in the current situation, we have an internal electrode, across which the normal component of the dielectric can (and will) jump. In NGSolve, we model this behavior by dividing the dielectric displacements into two parts, $d = d_1 + d_2$, where each part is defined in either the upper or the lower part. Additionally, we restrict the high-order shape functions to those that are divergence free, which allows to use only one degree of freedom per element for the electric potential
\begin{verbatim}
Sigma = HDivDiv(mesh, dirichlet=[2,3,5,6,7], order=k )
V = HCurl(mesh, dirichlet=[1], order=k)
Phi = L2(mesh, order=0 )
D1 = HDiv(mesh, order=k, dirichlet=[1,5,6,7], definedon=[1], hodivfree=True)
D2 = HDiv(mesh, order=k, dirichlet=[1,5,6,7], definedon=[2], hodivfree=True)
X = FESpace([Sigma, V, Phi,D1,D2] )

U = GridFunction(X)
Stress, Disp, Pot, DielD1, DielD2 = U.components
\end{verbatim}

The prescribed electric potential is now included into the right hand side of the variational equation \eqref{eq:V21}--\eqref{eq:divD2}, and no homogenization is needed.
Again, stiffness matrix and load vector are defined symbolically, reading
\begin{verbatim}
sigma,u,phi,d1,d2 = X.TrialFunction()
d_sigma,d_u,d_phi,d_d1,d_d2 = X.TestFunction()

a = BilinearForm(X, symmetric= False)
a += SymbolicBFI( InnerProduct (vec(sigma), SD*vec(d_sigma) ) )
a += SymbolicBFI( InnerProduct (g*vec(sigma), d_d1)\
                  +InnerProduct(g*vec(d_sigma), d1) )
a += SymbolicBFI( InnerProduct (g*vec(sigma), d_d2)\
                  +InnerProduct(g*vec(d_sigma), d2) )
a += SymbolicBFI( d1.Deriv()*d_phi + d_d1.Deriv()*phi )
a += SymbolicBFI( d2.Deriv()*d_phi + d_d2.Deriv()*phi )
a += SymbolicBFI( -InnerProduct( epsilonTinv*d1, d_d1) )
a += SymbolicBFI( -InnerProduct( epsilonTinv*d2, d_d2) )
a += SymbolicBFI( InnerProduct(div(sigma),d_u)+InnerProduct(div(d_sigma), u))
a += SymbolicBFI(-InnerProduct(sigma*n,tang(d_u))\
                 -InnerProduct(d_sigma*n,tang(u)), element_boundary=True)

phi_bd = CoefficientFunction([0,75,75,0,0,0,0])

f = LinearForm(X)
f += SymbolicLFI ( d_d1.Trace()*n*phi_bd, VOL_or_BND = BND )
f += SymbolicLFI ( d_d2.Trace()*n*phi_bd, VOL_or_BND = BND )
\end{verbatim}

Assembling and solving the linear system is the same as for the potential based method. The electric field is now evaluated using the material law.

\begin{verbatim}
E = -g*vec(Stress) + epsilonTinv*(DielD1 + DielD2)
\end{verbatim}

\rev{
\section{Conclusion}

In the present paper, we have introduced a non-standard finite element method for the simulation of piezoelectric materials under the assumptions of Voigt's linear theory. The method is a mixed method and includes stresses, namely the normal component of the stress vector, as independent unknowns. In a second variant, also the dielectric displacement is added as unknown field. In the purely elastic case, the elements have been shown to be locking-free when very flat, which makes them feasible for the discretization of flat piezoelectric structures. When dielectric displacements are discretized as well, the number of degrees of freedom can be reduced since only divergence-free shape functions need to be used. We present numerical results that indicate that also the flat piezoelectric elements converge at optimal order, as long as the electric potential is at least quadratic (in variant V1) or the dielectric displacement assumed linear (in variant V2). When using the elements in a homogenization procedure for a $d_{15}$ MFC, we see that good results are obtained for very coarse discretizations.  Due to the absence of locking, it is possible to discretize even the electrodes of thickness $\SI{2}{\micro\meter}$. All elements are available in the open-source finite element package Netgen/NGSolve. In the Appendix, an exemplary python script for the implementation of a bimorph beam is presented.

}
\section{Funding}
This work has been supported by the Linz Center of Mechatronics (LCM) in the framework of the Austrian COMET-K2 program.
Martin Meindlhumer acknowledges support of Johannes Kepler University Linz, Linz Institute of Technology (LIT).

\bibliographystyle{plain}
\bibliography{Smart} 


\end{document}